\documentclass[final,leqno]{siamltex}
\usepackage{amssymb}
\usepackage{amsmath}
\usepackage{color}
\usepackage{enumerate}
\usepackage{pgfplots}
\usepackage{tikz}
\usetikzlibrary{shapes.geometric}

\usepackage{graphicx}
\newcommand{\cT}{\mathcal{T}}

\newcommand{\cM}{\mathcal{M}}

\newcommand{\NJ}{\mathbb{N}_J^d}

\newcommand{\bx}{\mathbf{x}}

\newcommand{\del}{\delta}

\newcommand{\bh}{\mathbf{h}}

\newtheorem{remark}{Remark}[section]

\begin{document}

\title{Numerical moment stabilization of central difference approximations for linear stationary reaction-convection-diffusion equations 
with applications to stationary Hamilton-Jacobi equations }
 
\author{
Thomas Lewis\thanks{Department of Mathematics and Statistics, 
The University of North Carolina at Greensboro, 
Greensboro, NC 27402, U.S.A.  {\tt tllewis3@uncg.edu}.
The work of this author was partially supported by the NSF grant DMS-2111059.}
\and
Xiaohuan Xue\thanks{Department of Mathematics and Statistics, 
The University of North Carolina at Greensboro, 
Greensboro, NC 27402, U.S.A.  {\tt x\_xue2@uncg.edu}.
The work of this author was partially supported by the NSF grant DMS-2111059.}
}

\maketitle

\begin{abstract} 
Linear stationary reaction-convection-diffusion equations with Dirichlet boundary conditions 
are approximated using a simple finite difference method 
corresponding to central differences and the addition of a high-order stabilization term called a numerical moment.    
The focus is on convection-dominated equations, and the formulation for the method is motivated by  
various results for fully nonlinear problems.  The method features higher-order local truncation errors than monotone methods consistent with the 
use of the central difference approximation for the gradient. Stability and rates of convergence are derived in the $\ell^2$ norm for the constant-coefficient case.
Numerical tests are provided to compare the new methods to monotone methods.  
The methods are also tested for stationary Hamilton-Jacobi equations 
where they demonstrate higher rates of convergence than the Lax-Friedrich's method 
when the underlying viscosity solution is smooth 
and comparable performance when the underlying viscosity solution is not smooth.  
\end{abstract}



\pagestyle{myheadings}
\thispagestyle{plain}
\markboth{T. LEWIS and X. XUE}{Numerical moment stabilization of reaction-convection-diffusion equations}

%%%%%%%%%%%%%%%%%%%%%%%%%%%%%%%%%%%%%%%%%%%%%%%%%%
%%%%%%%%%%%%%%%%%%%%%%%%%%%%%%%%%%%%%%%%%%%%%%%%%%
%%%%%%%%%%%%%%%%%%%%%%%%%%%%%%%%%%%%%%%%%%%%%%%%%%
%%%%%%%%%%%%%%%%%%%%%%%%%%%%%%%%%%%%%%%%%%%%%%%%%%

\section{Introduction}

In this paper we explore the impact of adding a numerical moment stabilization term to a simple central-difference approximation scheme 
for stationary linear reaction-convection-diffusion equations of the form 
\begin{subequations} \label{RCD}
\begin{alignat}{2}
	L^{\epsilon}[u] \equiv -\epsilon \Delta u + \mathbf{b} \cdot \nabla u + c u & = f && \text{in } \Omega , \label{RCD_pde} \\ 
	u & = g \qquad && \text{on } \partial \Omega , \label{RCD_bc}
\end{alignat}
\end{subequations}
where 
$\epsilon$ is a positive constant; 
$\Delta u(\bx)$ is the Laplacian of $u$ at $\bx$; 
$\nabla u (\bx)$ is the gradient of $u$ at $\bx$;  
$\Omega \subset \mathbb{R}^{d}$ is a $d$-rectangle; 
$\mathbf{b} : \Omega \to \mathbb{R}^d$ is bounded; 
$c : \Omega \to \mathbb{R}$ is bounded and nonnegative; 
$f$ is bounded on $\Omega$; 
and $g$ is continuous on $\partial \Omega$.  
The problem is a singularly-perturbed  first-order problem when $\epsilon \ll 1$, i.e., the convection-dominated case.  
We will derive stability results and rates of convergence in $\ell^2$ in the special case when $\mathbf{b}$ is constant-valued.  
We will show that the results can be extended to the degenerate case $\epsilon = 0$ corresponding to $u$ being the viscosity solution of the reaction-convection equation 
$\mathbf{b} \cdot \nabla u + cu = f$ with the boundary condition $u=g$ satisfied in the viscosity-sense as long as there exists a constant 
$c_0$ such that $c \geq c_0 > 0$.  
By using central difference and a high-order stabilization technique, 
the methods can achieve higher-order local truncation errors than monotone methods 
while also avoiding the need for the mesh to resolve the diffusion constant 
when $\epsilon \ll 1$.  
The focus of this paper is to explore the impact of adding a numerical moment as a high-order stabilization term.  

The use of the central difference operator when discretizing $\nabla u$ 
instead of using a standard upwinding technique 
has two primary motivations.  
The first motivation is the fact that the matrix representation of the central difference approximation for $u_{x_i}$ 
is banded and skew-symmetric.  
Such a property can have major benefits when analyzing the numerical stability for the dynamic version of \eqref{RCD} 
with variable coefficients as seen in \cite{Hairer_Iserles_2016}.  
The central difference approximation of a partial derivative is the highest-order approximation that is skew-symmetric for 
uniform grids and Dirichlet boundary conditions.  
The second motivation is the second order accuracy of the central difference operator.  
Any monotone approximation of \eqref{RCD} is inherently limited to first-order accuracy 
in the limiting case $\epsilon \to 0^+$.
Unfortunately, it is not practical to directly use the central difference gradient approximation in \eqref{RCD} when $\epsilon \ll 1$.    
Suppose the mesh for $\Omega$ has uniform spacing $h_i$ in each direction $x_i$, and let $h = \max h_i$.  
Then, in order for the central difference operator to yield a monotone method when $c = 0$, 
the mesh must satisfy $h \leq 2 \epsilon / \| \mathbf{b} \|_{C^0(\Omega)}$ as seen in \cite{Birkhoff_Gartland_Lynch90}.  
Furthermore, in the degenerate case when $\epsilon = 0$, $\mathbf{b}$ is constant-valued (and nonzero), and $c = 0$, 
the matrix representation for \eqref{RCD} using central differences to approximate $\nabla u$ is singular when the number of rows is odd 
due to the fact that the matrix is skew-symmetric.  
Given the second-order accuracy of central difference, we investigate whether adding a simple stabilization term 
can allow for the practical use of the central difference operator when $\epsilon \ll 1$ 
as an alternative to using methods with reduced order such as upwinding schemes or 
more involved approximation techniques to formally achieve higher orders of accuracy.

When $\epsilon \ll 1$, equation \eqref{RCD} is closely related to the vanishing viscosity approximation of the stationary Hamilton-Jacobi (HJ) equation 
\begin{subequations} \label{HJ}
\begin{alignat}{2}
	H[u] \equiv H(\nabla u , u , \bx) & = 0 && \text{in } \Omega , \label{HJ_pde} \\ 
	u & = g \qquad && \text{on } \partial \Omega , \label{HJ_bc}
\end{alignat}
\end{subequations}
where $H$ is globally Lipschitz with respect to the arguments $\nabla u$ and $u$ 
and nondecreasing with respect to $u$.  
If $H$ also satisfies a comparison principle, then \eqref{HJ} has a unique viscosity solution in the space 
$C(\overline{\Omega})$. 
Furthermore, the solution $u^\epsilon$ of the perturbed problem 
\begin{subequations}\label{HJeps}
\begin{alignat}{2}
    H^\epsilon[u^\epsilon] \equiv - \epsilon \Delta u^\epsilon + H(\nabla u^\epsilon, u^\epsilon , \bx) & = 0 && \text{in } \Omega , \\ 
	u^\epsilon & = g \qquad && \text{on } \partial \Omega 
\end{alignat}
\end{subequations} 
converges to the solution $u$ of \eqref{HJ} in $L^\infty (\Omega)$ at a rate of $\mathcal{O}(\sqrt{\epsilon})$, 
where the Dirichlet boundary condition is understood in the viscosity sense (\cite{Crandall_Lions84}).  
Any linearization of \eqref{HJeps} would result in an equation of the form \eqref{RCD}.  
Thus, the ability to reliably approximate \eqref{RCD} with a higher-order method could have implications 
for approximating the viscosity solution of \eqref{HJ} with higher-order methods.  
We test our numerical moment stabilized versions of the central difference method for various HJ equations 
to further explore the practicality of the methods.  
Since the convection field is not known when approximating a solution to \eqref{HJ} or \eqref{HJeps}, 
the formulation of upwinding finite difference (FD) methods typical for approximating \eqref{RCD} is not possible 
with the added difficulty that the inflow and outflow boundaries depend on the function $u$ that is input into the operator $H[u]$.  
Thus, we explore a non-upwinding, high-order technique paired with a simple stabilizer to see how the 
boundary layer near an outflow boundary impacts the accuracy of the method throughout the interior of the domain.  

Convection-diffusion equations and stationary HJ equations arise from various scientific applications including 
optimal control, wave propagation, geometric optics, multiphase flow, image processing, etc. (cf. \cite{Osher, Sethian99} and the references therein). 
The numerical approximations to the solutions to time-dependent and stationary convection-diffusion equations and 
HJ equations have been significantly studied 
(cf. \cite{Barth, Brezzi_Russo, Crandall_Lions84, CockburnD, Hu, Huang, Kao, OsherShu, Sethian03, Shu, Tadmor, Tsai} and the references therein) 
as they are vital in understanding application problems. 
Monotone FD methods under the framework of Crandall and Lions in \cite{Crandall_Lions84} 
such as the Lax-Friedrich's methods borrowed from the approximation theory for nonlinear hyperbolic conservation laws 
(see \cite{Tadmor87}) have formed an analytic foundation for approximating viscosity solutions of HJ equations. 
Stationary Hamilton-Jacobi equations can be viewed  as degenerate fully nonlinear elliptic problems. 
The theory for monotone FD methods was further extended to elliptic problems by Barles and Souganidis in \cite{Barles}, 
which guarantees convergence to the underlying viscosity solution if the numerical scheme is 
monotone, admissible, consistent, and stable. 
Unfortunately, monotone methods for approximating \eqref{HJ} are limited to first-order accuracy due to the Godunov barrier (see \cite{Tadmor}). 
The main idea in this paper will be to treat \eqref{HJ} as a degenerate second order operator and  
introduce a {numerical moment} that will serve 
as a higher order stabilizer, as seen in \cite{FDhjb} and \cite{Kellie} for approximating fully nonlinear uniformly elliptic problems 
with non-monotone central difference approximations.  
Using a numerical moment will allow us to formally break the first-order accuracy barrier 
while also guaranteeing stability in certain cases.  
Note that several approaches have been used to design high-order methods to overcome the Godunov barrier 
for nonlinear problems 
using various types of spatial discretization techniques such as FD, finite elements, and discontinuous Galerkin methods 
(see \cite{Barth, Shu, centralDG_Li,Xu,Yan_Osher11} and the reference therein);  
however, convergence results typically require additional assumptions or regularity, are open due to the lack of monotonicity, 
or hold for dynamic problems but do not extend to stationary problems.  
The formulations in this paper are simple, and the numerical moment stabilization could easily be combined with 
more advanced techniques.  

This paper is written as a direct complement to the paper \cite{HJfiltered} which developed admissibility, stability, and convergence 
results for a related high-order FD method for approximating viscosity solutions of stationary 
Hamilton-Jacobi equations.  
In \cite{HJfiltered}, a cutoff operator was introduced to ensure admissibility and $\ell^\infty$-norm stability for a modified version of 
the scheme in this paper.  
By assuming constant-coefficients, we directly prove admissibility and $\ell^2$-norm stability for the unmodified, high-order scheme.  
The analytic techniques in this paper with an emphasis on $\ell^2$ results 
more closely resemble and generalize the techniques originally developed in \cite{FDhjb}.  
The series of papers shows the power of using a numerical moment as a stabilizer when approximating viscosity solutions 
of fully nonlinear boundary value problems 
while exploring various properties of the numerical moment itself.  
This paper also complements the paper \cite{DWDG_convection} which analyzes a novel discontinuous Galerkin method for approximating 
linear, convection-dominated problems.  
The linear case is well-studied with several references provided in \cite{DWDG_convection} as well as a more thorough discussion about \eqref{RCD} 
when the (potentially variable) convective velocity $\mathbf{b}$ is known so that an upwind direction is fixed.  
The goal in this paper is to test the impact of the numerical moment without taking advantage of the fixed convective velocity 
since it would be unknown when approximating \eqref{HJ} or \eqref{HJeps}.

The remainder of this paper is organized as follows. 
Some common notation, our mesh and difference operator notation, and background for the Lax-Friedrich's method and monotone methods 
for stationary HJ equations is provided in Section~\ref{preliminaries_sec}. 
The formulation of our stabilized central-difference FD methods for approximating \eqref{RCD}, \eqref{HJ}, and \eqref{HJeps} 
can be found in Section~\ref{formulation_sec}. 
Section~\ref{consistency_sec} will examine the consistency of the schemes 
and Section~\ref{moment_sec} will provide some additional properties of the numerical moment.   
Section~\ref{l2_sec} provides a complete $\ell^2$ analysis of the scheme for linear constant-coefficient problems.  
The schemes are numerically tested in Section~\ref{numerics_sec} for various choices of the numerical moment.  
Both linear and nonlinear problems are tested.  
In addition to exploring rates of convergence when approximating smooth PDE solutions, 
we also test the performance of the methods when the boundary condition leads to a discontinuity at the outflow boundary 
as $\epsilon \to 0^+$ or when the viscosity solution has a corner.  
Lastly, some concluding remarks and future directions are provided in Section~\ref{conc_sec}.  
We emphasize that this paper focuses primarily on the consistency and $\ell^2$  stability analysis for the linear problem with constant coefficients.  
When combined with the results in \cite{HJfiltered}, the papers provide a foundation from which we hope to develop a 
complete convergence analysis for the nonlinear problem.

%%%%%%%%%%%%%%%%%%%%%%%%%%%%%%%%%%%%%%%%%%%%%%%%%%
%%%%%%%%%%%%%%%%%%%%%%%%%%%%%%%%%%%%%%%%%%%%%%%%%%
%%%%%%%%%%%%%%%%%%%%%%%%%%%%%%%%%%%%%%%%%%%%%%%%%%
%%%%%%%%%%%%%%%%%%%%%%%%%%%%%%%%%%%%%%%%%%%%%%%%%%

\section{Preliminaries} \label{preliminaries_sec}

In this section we introduce the function space notation and difference operator notation that will be 
used throughout the paper.  
We also introduce the Lax-Friedrich's method as motivation for the numerical moment stabilization 
technique.  
The discussion will focus on monotone methods for approximating the stationary HJ equation \eqref{HJ}, 
but it can also be applied to upwind methods for approximating \eqref{RCD}.  
The background and conventions introduced in this section can also be found in \cite{HJfiltered}.  
  
%%%%%%%%%%%%%%%%%%%%%%%%%%%%%%%%%%%%%%%%%%%%%%%%%%

\subsection{Difference operators}

We introduce various difference operators for approximating first and second order partial derivatives 
using notation similar to \cite{FDhjb}.  
Assume $\Omega$ is a $d$-rectangle, i.e., 
$\Omega = \left( a_1 , b_1 \right) \times \left( a_2 , b_2 \right) \times \cdots \times 
		\left( a_d , b_d \right)$.    
We will only consider grids that are uniform in each coordinate $x_i$, $i = 1, 2, \ldots, d$.  
Let $J_i$ be a positive integer and $h_i = \frac{b_i-a_i}{J_i-1}$ for $i = 1, 2, \ldots, d$. 
Define $\mathbf{h} = \left( h_1, h_2, \ldots, h_d \right) \in \mathbb{R}^d$, 
$h = \max_{i=1,2,\ldots,d} h_i$, $h_* = \min_{i=1,2,\ldots,d} h_i$, $J = \prod_{i=1}^d J_i$, and   
$\NJ = \{ \alpha = (\alpha_1, \alpha_2, \ldots, \alpha_d) 
\mid 1 \leq \alpha_i \leq J_i, i = 1, 2, \ldots, d \}$.  Then, $\left| \NJ \right| = J$. 
We partition $\Omega$ into $\prod_{i=1}^d \left(J_i-1 \right)$ sub-$d$-rectangles with grid points
$\bx_{\alpha} = \Big(a_1+ (\alpha_1-1)h_1 , a_2 + (\alpha_2-1) h_2 , \ldots , 
		a_d + (\alpha_d-1) h_d \Big)$
for each multi-index $\alpha \in \NJ$.
We call $\cT_{\mathbf{h}}=\{ \bx_{\alpha} \}_{\alpha \in \NJ}$ 
a grid (set of nodes) for $\overline{\Omega}$. 

Let $\left\{ \mathbf{e}_i \right\}_{i=1}^d$ denote the canonical basis vectors for $\mathbb{R}^d$. 
We introduce an interior grid that removes the layer of grid points adjacent to the boundary of $\Omega$.  
To this end, let $\mathring{\cT_{\bh}} \subset \cT_{\bh}$ such that 
\[
	\mathring{\cT_{\bh}} \equiv 
	\left\{ \bx_\alpha \in \cT_{\bh} \cap \Omega \mid \bx_{\alpha \pm \mathbf{e}_i} \notin \partial \Omega 
	\text{ for all } i = 1,2,\ldots,d \right\} . 
\]
We let $\mathcal{S}_{h_i} \subset \cT_{\bh} \cap \partial \Omega$ denote the portion of the boundary 
normal to $\mathring{\cT_{\bh}}$ in the $\pm \mathbf{e}_i$ direction, 
i.e., 
\begin{align}\label{Sh_grid}
	\mathcal{S}_{h_i} \equiv \big\{ \bx_\alpha \in \cT_{\bh} \cap \partial \Omega \mid & \,  
		\bx_{\alpha} + h_i \mathbf{e}_i \in \cT_{\bh} \cap \Omega \text{ or } 
		\bx_{\alpha} - h_i \mathbf{e}_i \in \cT_{\bh} \cap \Omega \big\} 
\end{align}
for all $i \in \{1,2,\ldots,d\}$.  
Note that $\mathcal{S}_{h_i} \neq \emptyset$ for at most 
one $i \in \{1,2,\ldots,d\}$ at a given node $\bx_\alpha \in \cT_{\bh} \cap \partial \Omega$.  

Define the (first order) forward and backward difference operators by
\begin{equation*} %\label{fd_x}
	\del_{x_i,h_i}^+ v(\bx)\equiv \frac{v(\bx + h_i \mathbf{e}_i) - v(\bx)}{h_i},\qquad
	\del_{x_i,h_i}^- v(\bx)\equiv \frac{v(\bx)- v(\bx-h_i \mathbf{e}_i)}{h_i}
\end{equation*}
for a function $v$ defined on $\mathbb{R}^d$ and 
\[
	\del_{x_i,h_i}^+ V_\alpha \equiv \frac{V_{\alpha + \mathbf{e}_i} - V_\alpha}{h_i} , \qquad
	\del_{x_i,h_i}^- V_\alpha \equiv \frac{V_\alpha- V_{\alpha - \mathbf{e}_i}}{h_i}
\]
for a grid function $V$ defined on the grid $\mathcal{T}_{\mathbf{h}}$.  
We also define the following (second order) central difference operator: 
\begin{equation*} %\label{fd_xc}
	\delta_{x_i, h_i}  \equiv \frac{1}{2} \left( \delta_{x_i, h_i}^+ + \delta_{x_i, h_i}^- \right) 
         = \frac{v(\mathbf{x}+h_i \mathbf{e}_i) - v(\mathbf{x}-h_i \mathbf{e}_i)}{2 h_i}
\end{equation*}
and the ``sided" and central gradient operators $\nabla_{\mathbf{h}}^+$, 
$\nabla_{\mathbf{h}}^-$, and $\nabla_{\mathbf{h}}$ by 
\begin{align*} %\label{discrete_grad_def}
	\nabla_{\mathbf{h}}^\pm \equiv \bigl[ \delta^\pm_{x_1, h_1} , \delta_{x_2, h_2}^\pm , \cdots , 
		\delta_{x_d, h_d}^\pm \bigr]^T, \quad   %\label{sided_grad_def} \\ 
	\nabla_{\mathbf{h}} \equiv \bigl[ \delta_{x_1, h_1} , \delta_{x_2, h_2} , \cdots , 
		\delta_{x_d, h_d} \bigr]^T. %\label{central_grad_def} 
\end{align*}

Define the (second order) central difference operator for approximating second order 
partial derivatives by 
\begin{equation} \label{fd_xx}
	\del_{x_i,h_i}^2 v(\bx) 
	\equiv \frac{v(\bx + h_i \mathbf{e}_i) - 2 v(\bx) + v(\bx - h_i \mathbf{e}_i)}{h_i^2} 
\end{equation}
for a function $v$ defined on $\mathbb{R}^d$ and 
\[
	\del_{x_i,h_i}^2 V_\alpha 
	\equiv \frac{V_{\alpha + \mathbf{e}_i} - 2 V_\alpha + V_{\alpha - \mathbf{e}_i}}{h_i^2} 
\]
for a grid function $V$ defined on the grid $\mathcal{T}_{\mathbf{h}}$.  
We then define the (second order) central discrete Laplacian operator by 
\begin{equation*} %\label{LaplacianH}
	\Delta_{\mathbf{h}} \equiv \sum_{i=1}^d \delta_{x_i, h_i}^2 . 
\end{equation*}
In the formulation we will also consider the ``staggered" (second order) central difference operators 
$\delta_{x_i, 2h_i}^2$ defined by replacing $\mathbf{h}$ with $2\mathbf{h}$ in 
\eqref{fd_xx}.  
The ``staggered" operators have a 5-point stencil in each Cartesian direction and are defined using 
3 nodes.  
Note that ``ghost-values" may need to be 
introduced in order for the staggered difference operators to be well-defined at interior nodes adjacent to 
the boundary of $\Omega$.

%%%%%%%%%%%%%%%%%%%%%%%%%%%%%%%%%%%%%%%%%%%%%%%%%%

\subsection{The Lax-Friedrich's method and monotonicity}

In this section we introduce the convergence framework of Crandall and Lions and the idea of monotonicity 
for approximating solutions to \eqref{HJ}.  
We then introduce the Lax-Friedrich's method as an example of a 
method that falls within the Crandall and Lions framework.  
We will again consider the Lax-Friedrich's method in Section~\ref{numerics_sec} 
as a baseline for gauging the performance of our proposed method 
in a series of numerical tests for stationary HJ equations.  
A closer look at the Lax-Friedrich's method will help motivate our use of a higher-order stabilization technique 
for approximating \eqref{RCD} and \eqref{HJ}. 

The framework of Crandall and Lions ensures the convergence of methods for approximating 
viscosity solutions of \eqref{HJ} when using consistent and monotone methods.  
Suppose a scheme for approximating \eqref{HJ} has the form 
\[
	\widehat{H} \left( \nabla_\mathbf{h}^+ U_\alpha , \nabla_\mathbf{h}^- U_\alpha , U_\alpha , x_\alpha \right) = 0 , 
\]
where $\widehat{H}$ is called a numerical Hamiltonian.  
Assuming $\widehat{H}$ is continuous, we say $\widehat{H}$ is {\it consistent} if 
$\widehat{H} \left( \mathbf{q} , \mathbf{q} , v , \bx \right) = H(\mathbf{q}, v , \bx)$ 
for any vector $\mathbf{q} \in \mathbb{R}^d$, $v \in \mathbb{R}$, and $\bx \in \Omega$.  
For the stationary problem, 
we say $\widehat{H}$ is {\it monotone} on $[-R,R]$ if it is 
nondecreasing with respect to the node $U_\alpha$ and 
nonincreasing with respect to each 
node $U_{\alpha^\prime}$ such that $\bx_{\alpha^\prime} \neq \bx_\alpha$ is in the local stencil 
centered at $\bx_\alpha$ 
whenever $\left| \delta_{x_i, h_i}^+ U_{\alpha^\prime} \right| \leq R$.  
Consequently, the scheme is monotone if the numerical Hamiltonian $\widehat{H}$ 
is nonincreasing with respect to $\nabla_{\mathbf{h}}^+ U_\alpha$ 
and nondecreasing with respect to $\nabla_{\mathbf{h}}^- U_\alpha$ and $U_\alpha$.  

The Lax-Friedrich's numerical Hamiltonian is defined by 
\begin{align} \label{LF_scheme}
	\widehat{H}_{\text{LF}} U_\alpha 
	& \equiv H \left( \frac12 \nabla_{\mathbf{h}}^+ U_\alpha + \frac12 \nabla_{\mathbf{h}}^- U_\alpha , 
			U_\alpha , \bx_\alpha \right) 
		- \vec{\beta} \cdot \left( \nabla_{\mathbf{h}}^+ - \nabla_{\mathbf{h}}^- \right) U_\alpha \\ 
	\nonumber & = H \left( \nabla_{\mathbf{h}} U_\alpha , U_\alpha , \bx_\alpha \right) 
		- \vec{\beta} \cdot \left( \nabla_{\mathbf{h}}^+ - \nabla_{\mathbf{h}}^- \right) U_\alpha , 
\end{align}
where $\vec{\beta} \geq \mathbf{0}$ using the standard partial ordering for vectors.  
The numerical Hamiltonian is clearly consistent.  
Observe that, for $H$ Lipschitz and increasing with 
respect to $u$, 
if each component of $\vec{\beta}$ is sufficiently large, 
then the Lax-Friedrich's numerical Hamiltonian is monotone. 
The term $- \vec{\beta} \cdot \left( \nabla_{\mathbf{h}}^+ - \nabla_{\mathbf{h}}^- \right) U_\alpha$ is called a 
{\it numerical viscosity} due to the fact 
\[- \vec{\beta} \cdot \left( \nabla_{\mathbf{h}}^+ - \nabla_{\mathbf{h}}^- \right) U_\alpha 
	= - \sum_{i=1}^d \beta_i \frac{U_{\alpha + \mathbf{e}_i} - 2 U_\alpha + U_{\alpha - \mathbf{e}_i}}{h_i} 
	 = - \sum_{i=1}^d \beta_i h_i \delta_{x_i, h_i}^2 U_\alpha . 
\]
Choosing $\vec{\beta} = \beta \vec{1}$ for some constant $\beta \geq 0$, 
we have $\vec{\beta} \cdot \left( \nabla_{\mathbf{h}}^+ - \nabla_{\mathbf{h}}^- \right) U_\alpha$ 
is a (second-order) 
central difference approximation of $-\beta \Delta u(\bx_\alpha)$ 
scaled by $h$ when the mesh is uniform. 
We also have the method is inherently limited to first order accuracy due to the numerical viscosity.  
If $H$ is linear and we choose $\beta_i = \frac12 \left| \frac{\partial H}{\partial u_{x_i}} \right|$, 
then the particular Lax-Friedrich's method would correspond to the 
upwind method when approximating the linear problem \eqref{RCD} with $\epsilon = 0$.  

\begin{remark}
Two aspects of the Lax-Friedrich's method will serve as motivation for using a numerical moment stabilizer.  
First, the Lax-Friedrich's scheme is analogous to using a central difference approximation of the gradient and choosing $\epsilon = \beta h$ 
in \eqref{HJeps} as a stabilizer.  
The numerical moment will be used in lieu of a numerical diffusion to achieve higher orders of accuracy.  
Second, the method in the context of the Crandall and Lions framework naturally incorporates two discrete gradient operators.  
The proposed method in Section~\ref{formulation_sec} will naturally incorporate two second order operators to form a high-order stabilizer.  
\end{remark}

We end this section by noting that the first order accuracy bound is a consequence of the monotonicity.  
Every convergent monotone FD scheme for \eqref{HJ} implicitly approximates the differential equation 
\begin{equation} \label{tadmor_eqn}
	- \beta h ``\Delta u" + H(\nabla u , u ,\bx) = 0 
\end{equation}
for sufficiently large and possibly nonlinear $\beta > 0$, 
where $- \beta h ``\Delta u"$ corresponds to a numerical viscosity (c.f. \cite{Tadmor}).  
Thus, we cannot increase the power of the coefficient $h$ for the numerical viscosity 
and maintain a monotone scheme.  
Similar observations hold when using the central difference approximation of $\nabla u$ in the 
linear reaction-convection-diffusion equation \eqref{RCD} with $\epsilon \ll 1$ unless a 
mesh condition is introduced for which $h = \mathcal{O}(\epsilon)$.  
To achieve a higher order scheme and to potentially avoid mesh conditions, 
we will necessarily abandon the monotonicity condition when introducing an alternative high-order stabilizer.

%%%%%%%%%%%%%%%%%%%%%%%%%%%%%%%%%%%%%%%%%%%%%%%%%%
%%%%%%%%%%%%%%%%%%%%%%%%%%%%%%%%%%%%%%%%%%%%%%%%%%
%%%%%%%%%%%%%%%%%%%%%%%%%%%%%%%%%%%%%%%%%%%%%%%%%%
%%%%%%%%%%%%%%%%%%%%%%%%%%%%%%%%%%%%%%%%%%%%%%%%%%

\section{Formulation} \label{formulation_sec}

We approximate the solution of \eqref{RCD}, \eqref{HJ}, or \eqref{HJeps} using a FD method with 
up to second order local truncation error.  
We will introduce two different choices for an auxiliary boundary condition that ensure the 
corresponding system of equations is well-defined.  
The main idea when formulating the new method 
will be to add a high-order stabilization term (called a numerical moment) that will formally 
allow for the choice of $h_*^2$ in \eqref{tadmor_eqn}. 
Such a choice will sacrifice the monotonicity assumption, but it is consistent 
with the vanishing viscosity interpretation in \eqref{HJeps}.    

The new method is defined as follows. 
Let $\gamma, \sigma \geq 0$ be constants and $p \in [0,1]$.  
The proposed FD method is defined as finding a grid function $U_\alpha$ such that 
\begin{subequations} \label{RCD_fd}
\begin{alignat}{2} 
	\widehat{L}_{\mathbf{h}}^\epsilon [ U_\alpha ] & = f(\bx_\alpha) 
		&& \text{if } \bx_\alpha \in \cT_{\bh} \cap \Omega , \label{RCD_fd1} \\ 
	U_\alpha & = g(\bx_\alpha) \quad && \text{if } \bx_\alpha \in \cT_{\bh} \cap \partial \Omega \label{RCD_fd2}  , \\ 
	B_{\bh} U_\alpha & = 0 \quad && \text{if } \bx_\alpha \in  \mathcal{S}_{h_i} \subset \mathcal{T}_{\mathbf{h}} \cap \partial\Omega \label{RCD_fd3} 
\end{alignat} 
\end{subequations} 
when approximating \eqref{RCD} 
or 
\begin{subequations} \label{HJeps_fd}
\begin{alignat}{2} 
	\widehat{H}_{\mathbf{h}}^\epsilon [ U_\alpha ] & = 0 
		&& \text{if } \bx_\alpha \in \cT_{\bh} \cap \Omega , \label{HJ_fd1} \\ 
	U_\alpha & = g(\bx_\alpha) \quad && \text{if } \bx_\alpha \in \cT_{\bh} \cap \partial \Omega \label{HJ_fd2} , \\ 
	B_{\bh} U_\alpha & = 0 \quad && \text{if } \bx_\alpha \in  \mathcal{S}_{h_i} \subset \mathcal{T}_{\mathbf{h}} \cap \partial\Omega \label{HJ_fd3} 
\end{alignat} 
\end{subequations} 
for all $\epsilon \geq 0$ when approximating \eqref{HJ} or \eqref{HJeps}, 
where the numerical operators $\widehat{L}_{\bh}^\epsilon$ and $\widehat{H}_{\bh}^\epsilon$ are defined by 
\begin{align*}
	\widehat{L}_{\mathbf{h}}^\epsilon [U_\alpha] 
	& \equiv - \sigma \sum_{i=1}^d h_i^2 \delta_{x_i, h_i}^2 U_\alpha - \epsilon \Delta_{\bh} U_\alpha
		+ \mathbf{b}(\bx_\alpha) \cdot \nabla_{\bh} U_\alpha + c(\bx_\alpha) U_\alpha 
		+ \gamma M_{\bh}^p U_\alpha , \\ 
	\widehat{H}_{\mathbf{h}}^\epsilon [U_\alpha] 
	& \equiv - \sigma \sum_{i=1}^d h_i^2 \delta_{x_i, h_i}^2 U_\alpha - \epsilon \Delta_{\bh} U_\alpha
		+ H \left( \nabla_{\bh} U_\alpha , U_\alpha , \bx_\alpha \right) 
		+ \gamma M_{\bh}^p U_\alpha , 
\end{align*} 
the numerical moment operator $M_{\bh}^p$ is defined by  
\[
	M_{\bh}^p U_\alpha \equiv \sum_{i=1}^d h_i^p \left( \delta_{x_i, 2h_i}^2 - \delta_{x_i, h_i}^2 \right) U_\alpha , 
\]
and the auxiliary boundary condition operator $B_{\bh}$ is defined by either 
\begin{align} \label{bc2a}
	B_{\bh} U_\alpha \equiv -\delta_{x_i, h_i}^2 U_\alpha  
\end{align} 
or 
\begin{align} \label{bc2b}
	 B_{\bh} U_\alpha \equiv -\delta_{x_i, h_i}^2 U_\alpha +  \delta_{x_i, h_i}^2 U_{\alpha'} , 
\end{align} 
where $\bx_{\alpha'} \in \{ \bx_{\alpha} - h_i \mathbf{e}_i , \bx_{\alpha} + h_i \mathbf{e}_i \}$ 
such that $\bx_{\alpha'} \in \cT_{\bh} \cap \Omega$.  
The auxiliary boundary condition accounts for the ghost values needed to evaluate $\delta_{x_i, 2h_i}^2 U_\alpha$ 
for nodes adjacent to the boundary in the $x_i$ direction.  
In general, the choice $p=1$ requires $\gamma$ to be chosen sufficiently large as will be characterized in the stability analysis 
and as seen for the modified scheme in \cite{HJfiltered}.  

We refer to the term $M_{\bh}^p U_\alpha$ as a numerical moment 
due to the fact that 
\begin{align} \label{moment_expansion}
	M_{\bh}^p U_\alpha 
	& = \sum_{i=1}^d \frac{U_{\alpha+2\mathbf{e}_i} - 4 U_{\alpha+\mathbf{e}_i} + 6 U_\alpha 
		- 4 U_{\alpha - \mathbf{e}_i} + U_{\alpha - 2 \mathbf{e}_i}}{4 h_i^{2-p}} \\ 
	\nonumber & = \frac14 \sum_{i=1}^d h_i^{2+p} \delta_{x_i, h_i}^2 \delta_{x_i, h_i}^2 U_\alpha  
\end{align}
implying $M_{\bh}^p U_\alpha$ is a (second-order) 
central difference approximation of the fourth order differential operator 
$\sum_{i=1}^d u_{x_i x_i x_i x_i}(\bx_\alpha)$ with each term scaled by  
constants proportional to $h^{2+p}$ when the mesh is quasi-uniform.  
As such, the proposed scheme is a direct realization of both the vanishing viscosity method (when $\sigma > 0$ and $\epsilon = 0$) 
and the vanishing moment method (see \cite{Feng_Neilan11}).  
By adding a vanishing moment term to the Lax-Friedrich's scheme, 
the coefficient on the vanishing viscosity term 
can been scaled by an additional power of $h$.  
The vanishing moment will be essential to deriving stability bounds that hold when $0 \leq \epsilon \ll 1$.  

A more general form of the numerical moment can be found in \cite{Kellie}.  
We could also use the auxiliary boundary $\Delta_{\bh} U_\alpha = 0$ on $\cT_{\bh} \cap \partial \Omega$ 
to be consistent with the choice for the auxiliary boundary condition in \cite{FDhjb} which also considered 
mixed partial derivatives.  
Since the operators $L^\epsilon$ and $H^\epsilon$ do not include mixed partial derivatives, we instead use the 
boundary conditions that only look at nonmixed partials in a normal direction.  
Formally, the auxiliary boundary condition extends the solution to the PDE along 
vectors normal to the boundary of $\Omega$ using either a linear or quadratic extension.  
We lastly note that the numerical moment stabilization in \cite{FDhjb} and \cite{Kellie} corresponded to the choice $p=0$ 
when approximating second order fully nonlinear problems.  
In this paper we also consider the rescaled numerical moment corresponding to the choice $p \in [0,1]$.  

%%%%%%%%%%%%%%%%%%%%%%%%%%%%%%%%%%%%%%%%%%%%%%%%%%
%%%%%%%%%%%%%%%%%%%%%%%%%%%%%%%%%%%%%%%%%%%%%%%%%%
%%%%%%%%%%%%%%%%%%%%%%%%%%%%%%%%%%%%%%%%%%%%%%%%%%
%%%%%%%%%%%%%%%%%%%%%%%%%%%%%%%%%%%%%%%%%%%%%%%%%%

\section{Consistency results} \label{consistency_sec}

In this section we consider the consistency of the proposed scheme \eqref{HJeps_fd} 
paired with the auxiliary boundary condition \eqref{bc2a} or \eqref{bc2b}.  
The following will only consider \eqref{HJeps_fd} since \eqref{RCD_fd} can be considered a special case when $H$ is linear.  
First, suppose that $\delta_{x_i, h_i}^2 U_\alpha \to p_i$, $\delta_{x_i, 2h_i}^2 U_\alpha \to p_i$, 
$\nabla_{\bh}^\pm U_\alpha \to \mathbf{v}$, 
$U_\alpha \to v$, and $\bx_\alpha \to \bx_0$ as $h \to 0$.  
Then, we have 
\[
	\widehat{H}_{\bh}^\epsilon [U_\alpha] \to - \epsilon \sum_{i=1}^d p_i + H(\mathbf{v}, v, \bx_0) , 
\]
and it follows that the method is consistent using the definition in \cite{Kellie} for fully nonlinear elliptic problems 
(which naturally extends the notions of consistency in \cite{Crandall_Lions84} and \cite{Barles} for the case when 
multiple second derivative approximations are used).  
We next analyze the local truncation errors as well as the global $\ell^2$ truncation error for the scheme \eqref{HJeps_fd} 
paired with \eqref{bc2a} or \eqref{bc2b}.  

Let $v \in C^5(\overline{\Omega})$, and choose $\bx_\alpha \in \cT_{\bh}$.  
Clearly the scheme has zero local truncation error if $\bx_\alpha \in \partial \Omega$ by the Dirichlet boundary condition
(assuming $g(\bx_\alpha) = v(\bx_\alpha)$).  
Thus, we only consider $\bx_\alpha \in \Omega$.  
In \cite{HJfiltered}, it is shown that the auxiliary boundary conditions create a boundary layer for the local truncation error 
of the scheme \eqref{HJeps_fd} 
when evaluating the numerical moment.  In particular, it is shown that 
\[
	\widehat{H}_{\mathbf{h}}^\epsilon[ v(\bx_\alpha)] 
	= -\epsilon \Delta v(\bx_\alpha) + H \left( \nabla v(\bx_\alpha) , v(\bx_\alpha) , \bx_\alpha \right) + \mathcal{O}(h^2) 
\]
for all $\bx_\alpha \in \mathring{\cT_{\bh}}$ and   
\[
	\widehat{H}_{\mathbf{h}}^\epsilon[ v(\bx_\alpha)] 
	= -\epsilon \Delta v(\bx_\alpha) + H \left( \nabla v(\bx_\alpha) , v(\bx_\alpha) , \bx_\alpha \right) + \mathcal{O}(h^q) 
\]
for all $\bx_\alpha \in (\cT_{\bh} \setminus \mathring{\cT}_{\bh} ) \cap \Omega$, 
where $q = p$ if \eqref{bc2a} is used to define $B_{\bh}$ 
and $q=p+1$ if \eqref{bc2b} is used to define $B_{\bh}$.  

We see that the local truncation error is largest for grid points adjacent to the boundary since $0 \leq p \leq 1$.  
We also see that the second auxiliary boundary condition \eqref{bc2b} yields a higher-order local truncation error that \eqref{bc2a}.  
A simple calculation reveals that the boundary condition \eqref{bc2b} is equivalent to assuming 
$u_{x_i x_i x_i} = 0$ using a (forward or backwards) first-order approximation of the third derivative enforced at the ghost point.  
Since the auxiliary boundary condition \eqref{bc2a} assumes the second derivative is zero-valued when crossing the boundary 
while the auxiliary boundary condition \eqref{bc2b} assumes the second derivative is constant through the boundary,  
we have \eqref{bc2a} would be exact if $v$ is linear while \eqref{bc2b} 
would be exact if $v$ is quadratic.  

Combining the above observations, we have the following results that are proved in \cite{HJfiltered}.  

\begin{lemma}\label{local_truncation_lemma}
The scheme \eqref{HJeps_fd} has a second order $\ell^\infty$ local truncation error over $\mathring{\cT_{\bh}}$.  
The boundary condition \eqref{bc2a} gives a $p$ order $\ell^\infty$ local truncation error over 
$( \cT_{\bh} \setminus \mathring{\cT_{\bh}} ) \cap \Omega$, 
and the boundary condition \eqref{bc2b} gives a $p+1$ order $\ell^\infty$ local truncation error over 
$( \cT_{\bh} \setminus \mathring{\cT_{\bh}} ) \cap \Omega$, 
where $p \in [0,1]$.  
\end{lemma}

\begin{corollary}
The numerical moment is exact for linear functions when using the auxiliary boundary condition \eqref{bc2a}, 
and it is exact for quadratic functions when using the auxiliary boundary condition \eqref{bc2b}.  
The numerical viscosity scaled by $h_i^2$ presents an $\mathcal{O}(h^2)$ consistency error, 
and it is only exact for linear functions.  
The scheme is exact for linear functions for either choice of auxiliary boundary condition.  
\end{corollary}

We lastly bound the global truncation error measured in a weighted $\ell^2$ norm 
given by $(\prod_{i=1,2,\ldots,d} h_i^{1/2} ) \| V \|_{\ell^2(\cT_{\bh} \cap \Omega)}$ for a grid function $V$.  
Thus, the weighted norm approximates the $L^2$-norm of the underlying function $v$.   
The global truncation error weights the single boundary layer against the interior second-order error.  
The numerical tests in Section~\ref{numerics_sec} produce exact errors that are consistent with or better than 
the global truncation error.  

\begin{theorem} \label{GTE_theorem}
The global truncation error of scheme \eqref{HJeps_fd} with boundary condition \eqref{bc2a} has order $p+\frac12$. 
The global truncation error of scheme \eqref{HJeps_fd} with boundary condition \eqref{bc2b} has order $\min\{2,p+\frac{3}{2}\}$.
\end{theorem}

\begin{proof}
Suppose $v \in C^5(\overline{\Omega})$.  
Then, using Lemma~\ref{local_truncation_lemma}
and the facts that there exists a constant $C_{\Omega}$ such that $\left| \mathring{\cT}_{\bh} \right| \leq C_\Omega \frac{1}{\prod_{i=1,2,\ldots,d} h_i}$ 
and $\left| (\cT_{\bh}\setminus \mathring{\cT_{\bh}}) \cap \Omega \right| \leq C_\Omega \frac{\sum_{i=1,2,\ldots,d}h_i}{\prod_{i=1,2,\ldots,d}h_i}$, 
%$\| 1_{\bh} \|_{\ell^2(\mathring{\cT}_{\bh})}^2 = \sum_{\bx_\alpha \in \mathring{\cT}_{\bh}} 1 \leq \frac{1}{\prod_{i=1,2,\ldots,d} h_i}$, 
there holds 
\begin{align*}
   &\left( \prod_{i=1,2,\ldots,d} h_i \right) \left\| H[v (\bx_\alpha)]- \widehat{H}_{\mathbf{h}} [v (\bx_\alpha)] \right\|^2_{\ell^2(\cT_{\bh} \cap \Omega)} \\
   & \qquad = \left( \prod_{i=1,2,\ldots,d} h_i \right) 
   	\left\| H[v(\bx_\alpha)] - \widehat{H}_{\mathbf{h}} [v (\bx_\alpha)] \right\|^2_{\ell^2(\mathring{\cT_{\bh}} \cap \Omega)} \\ 
	& \qquad \qquad + \left( \prod_{i=1,2,\ldots,d} h_i \right) 
		\left\| H[v(\bx_\alpha)]- \widehat{H}_{\mathbf{h}} [v (\bx_\alpha)] \right\|^2_{\ell^2((\cT_{\bh}\setminus \mathring{\cT_{\bh}}) \cap \Omega)} \\
   & \qquad \leq \left( \prod_{i=1,2,\ldots,d} h_i \right) \left[ C \frac{1}{ \prod_{i=1,2,\ldots,d} h_i} h^{4}
   	+ C \frac{\sum_{i=1,2,\ldots,d}h_i}{\prod_{i=1,2,\ldots,d}h_i} h^{2(p+r)}\right] \\
   &\qquad  = Ch^{4} + C \left(\sum_{i=1,2,\ldots,d}h_i\right) h^{2q}\leq Ch^{4} + Cdh^{2q+1} \\ 
   & \qquad = \mathcal{O} \left(h^{\min\{4,2q+1\}} \right)  
\end{align*}   
for some constant $C$ independent of $h$, $q=p$ for auxiliary boundary condition \eqref{bc2a}, 
and $q=p+1$ for auxiliary boundary condition \eqref{bc2b}.  
The result follows by taking the square root of both sides.  
\hfill
\end{proof}

\begin{remark}
For comparison, any monotone method would have a local truncation error of order 1 and a 
global truncation error of order 1.  
In Section~\ref{numerics_sec}, we will see that the proposed method often converges faster than the order of the global truncation error 
with observed rates greater than 1 even when using the boundary condition \eqref{bc2a} with $p=0$.  
\end{remark}

%%%%%%%%%%%%%%%%%%%%%%%%%%%%%%%%%%%%%%%%%%%%%%%%%%
%%%%%%%%%%%%%%%%%%%%%%%%%%%%%%%%%%%%%%%%%%%%%%%%%%
%%%%%%%%%%%%%%%%%%%%%%%%%%%%%%%%%%%%%%%%%%%%%%%%%%
%%%%%%%%%%%%%%%%%%%%%%%%%%%%%%%%%%%%%%%%%%%%%%%%%%

\section{Matrix properties of the numerical moment} \label{moment_sec}

In this section we explore various properties of the numerical moment.  
In particular, we explore monotonicity properties and spectral properties.  
We see that the numerical moment is a higher-order stabilization operator 
that preserves select matrix properties of the numerical viscosity. 

Choose $i \in \{1,2,\ldots,d \}$.  
Let $D_{i,0}^2$ denote a matrix representation of the discrete operator $-\delta_{x_i, h_i}^2$ with zero Dirichlet boundary condition 
over $\mathcal{S}_{h_i}$
and $D_{i,1}^2$ denote a matrix representation of the discrete operator $-\delta_{x_i, h_i}^2$ with zero Neumann boundary condition 
over $\mathcal{S}_{h_i}$
imposed using either forward or backward difference quotients to avoid the need for a ghost point 
(using the same ordering of the unknowns for both matrices).  
Then $D_{i,0}^2$ a symmetric positive definite M-matrix 
and $D_{i,1}^2$ is a symmetric nonnegative definite M-matrix.  
Furthermore, the matrix representation of $M_{\bh}^p$ corresponds to 
the matrix $\frac14 h_i^{p+2} D_{i,0}^2 D_{i,0}^2$ when using the auxiliary boundary condition \eqref{bc2a} 
and the matrix $\frac14 h_i^{p+2} D_{i,1}^2 D_{i,0}^2$ when using the auxiliary boundary condition \eqref{bc2b}. 
We immediately have that the matrix representation of $M_{\bh}^p$ with the auxiliary boundary condition \eqref{bc2a} 
is a monotone matrix that is symmetric positive definite.  

We now consider the matrix representation of $M_{\bh}^p$ with the auxiliary boundary condition \eqref{bc2b} 
more closely.  
The matrix itself is singular with a nullspace corresponding to $-\delta_{x_i, h_i}^2 V_\alpha$ being constant-valued for some grid function $V_\alpha$.  
However, for $\sigma > 0$, we have the stabilization term associated with the $x_i$ direction is 
$- \sigma h_i^2 \delta_{x_i, h_i}^2 U_\alpha + \frac14 \gamma h_i^{2+p} (-\delta_{x_i, h_i}^2) (-\delta_{x_i, h_i}^2) U_\alpha$ 
due to the numerical viscosity.  
Thus, we have the matrix 
$\left(\sigma h_i^2 I + \frac{\gamma}{4} h_i^{p+2} D_{i,1}^2 \right) D_{i,0}^2$ 
is associated with the numerical viscosity and numerical moment.  
The matrix is monotone since the matrix $\sigma h_i^2 I + \frac{\gamma}{4} h_i^{p+2} D_{i,1}^2$ is a nonsingular M-matrix.  
However, the matrix is no longer symmetric positive definite since the product does not commute.

%%%%%%%%%%%%%%%%%%%%%%%%%%%%%%%%%%%%%%%%%%%%%%%%%%
%%%%%%%%%%%%%%%%%%%%%%%%%%%%%%%%%%%%%%%%%%%%%%%%%%
%%%%%%%%%%%%%%%%%%%%%%%%%%%%%%%%%%%%%%%%%%%%%%%%%%
%%%%%%%%%%%%%%%%%%%%%%%%%%%%%%%%%%%%%%%%%%%%%%%%%%

\section{Admissibility, stability, and convergence for constant coefficients} \label{l2_sec}

In this section we analyze the scheme \eqref{RCD_fd} with auxiliary boundary condition \eqref{bc2a} in the special 
case when $\mathbf{b}$ is constant-valued.  
Note that the analysis technique 
does not extend to the choice of auxiliary boundary condition \eqref{bc2b} due to the fact the 
matrix corresponding to the numerical moment $M_{\bh}^p$ is not symmetric 
and its symmetrization is not symmetric nonnegative definite (as can be verified computationally).  
The analysis similarly would not extend to the more general case when $\mathbf{b}$ is not constant-valued 
or the nonlinear problem \eqref{HJeps_fd} 
unless $H$ has the special form $H(\nabla u , u , \bx) = \mathbf{b} \cdot \nabla u + r(u,\mathbf{x}) - f(\mathbf{x})$ 
for $r$ globally Lipschitz and nondecreasing with respect to $u$.  
However, there are cases for dynamic problems where stability results can be obtained when $\mathbf{b}$ is Lipschitz 
as seen in \cite{Hairer_Iserles_2016}.  

Let $\cT_{\bh}'$ denote the extended mesh including the ghost points 
and $J' = \left| \cT_{\bh}' \right|$.  
We define a mapping $\cM_\rho : \mathbb{R}^{J'} \to \mathbb{R}^{J'}$ 
that will be used to define a fixed-point iteration for solving 
\eqref{RCD_fd} and \eqref{bc2a} when $\mathbf{b}$ is constant-valued.  
We show the mapping is contractive in the $\ell^2$-norm for $\rho > 0$ sufficiently small. 
From there, we can derive an $\ell^2$-norm stability estimate for the solution to the proposed FD scheme 
as well as rates of convergence.   

Let $S(\cT_{\bh}')$ denote the space of grid functions defined over $\cT_{\bh}'$.
We define the mapping $\cM_\rho : S(\mathcal{T}_{\mathbf{h}}') \to S(\cT_{\bh}')$ by 
\begin{equation}\label{M_rho}
	\widehat{U}_\alpha 
	\equiv \cM_\rho U_\alpha , 
\end{equation}
where 
\begin{subequations} 
\begin{alignat}{2} 
	\widehat{U}_\alpha 
	& = U_\alpha 
		- \rho \widehat{L}_{\bh}^\epsilon [U_\alpha ] - \rho f(\bx_\alpha) \qquad 
		&& \text{if } \bx_\alpha \in \mathcal{T}_{\mathbf{h}} \cap \Omega , \\ 
	\widehat{U}_\alpha & = g(\bx_\alpha) \quad && \text{if } \bx_\alpha \in \cT_{\bh} \cap \partial \Omega , \label{Wbc1} \\ 
	\delta_{x_i, h_i}^2 \widehat{U}_\alpha & = 0 \quad && \text{if } \bx_\alpha \in \mathcal{S}_{h_i} . \label{Wbc2}
\end{alignat} 
\end{subequations} 
Clearly a fixed point of \eqref{M_rho} is a solution to the FD scheme \eqref{RCD_fd} and \eqref{bc2a} and vice versa.  
We show \eqref{M_rho} is a contraction for all $\rho > 0$ sufficiently small.  

\begin{lemma} \label{contraction_uh} 
Choose $U,V \in S(\cT_{\bh}')$, and let $\widehat{U} = \cM_\rho U$ and $\widehat{V} = \cM_\rho V$ 
for $\cM_\rho$ defined by \eqref{M_rho}.  
Then, there holds 
\[
	\| \widehat{U} - \widehat{V} \|_{\ell^2(\cT_{\bh})} \leq \left(1-\rho \frac{\epsilon \lambda_0 + c_0 + h_*^{2+p} \gamma \lambda_*}{2}\right) \| U - V \|_{\ell^2(\cT_{\bh})}
\] 
for all $\rho > 0$ sufficiently small, 
where $\lambda_0$ is the minimal eigenvalue of $-\Delta_{\bh}$, $c_0 = \inf_{\Omega} c$, and $\lambda_*$ is the minimal eigenvalue 
of $\frac14 \delta_{x_i, h_i}^2 \delta_{x_i, h_i}^2$ with Dirichlet boundary data enforced for the operator $\delta_{x_i, h_i}^2$ 
over $\mathcal{S}_{h_i}$ for $h_i = h_*$.  
\end{lemma}

\begin{proof}
Let $W = V - U$ and $\widehat{W} = \widehat{V} - \widehat{U}$.  
Then 
\begin{subequations} \label{V_U}
\begin{alignat}{2} 
	\widehat{W}_\alpha 
	& = (I - \rho \widehat{L}_{\bh}^\epsilon) [W_\alpha ] \qquad 
		&& \text{if } \bx_\alpha \in \mathcal{T}_{\mathbf{h}} \cap \Omega , \\ 
	\widehat{W}_\alpha & = 0 \quad && \text{if } \bx_\alpha \in \cT_{\bh} \cap \partial \Omega ,  \\ 
	\delta_{x_i, h_i}^2 \widehat{W}_\alpha & = 0 \quad && \text{if } \bx_\alpha \in \mathcal{S}_{h_i} 
\end{alignat} 
\end{subequations} 
by the linearity of $\widehat{L}_{\bh}^\epsilon$ and the boundary conditions.  

Let $J_0 = \left| \cT_{\bh} \cap \Omega \right|$ and 
$\widehat{\mathbf{W}}, \mathbf{W} \in \mathbb{R}^{J_0}$ denote vectorizations 
of $\widehat{W}$ and $W$, respectively.  
Define the matrices $C, L, M \in \mathbb{R}^{J_0 \times J_0}$ 
such that $C$ is a diagonal matrix corresponding to the values of $c(\bx_\alpha)$, 
$L$ corresponds to $-\Delta_{\bh}$ with Dirichlet boundary condition, 
and $M$ corresponds to $M_{\bh}^p$ with Dirichlet boundary condition and the auxiliary boundary condition, respectively, 
Notationally, we write $A \leq B$ for symmetric matrices $A, B$ if $B-A$ is symmetric nonnegative definite.  
Then, the matrices $C, L, M$ are symmetric nonnegative definite with 
\[
	C \geq c_0 I, \qquad L \geq \lambda_0 I , \qquad M \geq h_*^{2+p} \lambda_* I .  
\]
Lastly, define the symmetric positive definite matrices $D_{ii,0}$ corresponding to $-\delta_{x_i, h_i}^2$ 
with Dirichlet boundary conditions and the 
antisymmetric matrices $D_i$ corresponding to $\delta_{x_i, h_i}$ 
with Dirichlet boundary conditions.  

Using the above vectorizations, we have \eqref{V_U} is equivalent to 
\begin{align*}
\widehat{\mathbf{W}} & = \mathbf{W} - \rho \epsilon L \mathbf{W} - \rho C \mathbf{W} 
	- \rho \sigma \sum_{i=1}^d h_i^2 D_{ii,0} \mathbf{W} - \rho \gamma M \mathbf{W} 
	- \rho \sum_{i=1}^d b_i D_i \mathbf{W} \\ 
	& \equiv (I - \rho G_s - \rho G_a) \mathbf{W} 
\end{align*}
for 
$G_s \equiv \epsilon L + C + \gamma M + \rho \sigma \sum_{i=1}^d h_i^2 D_{ii,0}$ symmetric nonnegative definite 
and $G_a \equiv \sum_{i=1}^d b_i D_i$ antisymmetric since $b_i$ is a scalar and $D_i$ is antisymmetric for each $i$.  
There exists a constant $R_1$ such that  
\[
	\mathbf{0} \leq \frac12 I - \rho G_s \leq \left( \frac12 - \rho \epsilon \lambda_0 - \rho c_0 - \rho h_*^{2+p} \gamma \lambda_* \right) I 
\]
for all $0 < \rho < R_1$.  
Observe that 
\[
	(I - \rho G_a)^T (I - \rho G_a) = I - \rho G_a^T - \rho G_a + \rho^2 G_a^T G_a = I + \rho^2 G_a^T G_a 
\]
by the antisymmetry of $G_a$.   
Let $\kappa$ be the maximum eigenvalue of $G_a^T G_a$.  
Then, $\| I - \rho G_a \|_2 \leq \sqrt{ 1 + \rho^2 \kappa }$, 
and we have 
\begin{align*}
	\| \widehat{\mathbf{W}} \|_2 
	& \leq \| I - \rho G_s - \rho G_a \|_2 \| \mathbf{W} \|_2 \\ 
	& \leq \| \frac12 I - \rho G_s \|_2 \| \mathbf{W} \|_2 + \frac12 \left\| I - 2 \rho G_a \right\|_2 \| \mathbf{W} \|_2 \\ 
	& \leq \left( \frac12 - \rho \epsilon \lambda_0 - \rho c_0 - \rho h_*^{2+p} \gamma \lambda_* \right) \| \mathbf{W} \|_2 
		+ \frac12 \left\| I - 2 \rho G_a \right\|_2 \| \mathbf{W} \|_2 \\ 
	& \leq \left( \frac12 - \rho \epsilon \lambda_0 - \rho c_0 - \rho h_*^{2+p} \gamma \lambda_* + \frac12 \sqrt{1 + 4 \rho^2 \kappa} \right) \| \mathbf{W} \|_2 \\ 
	& = \left( \frac12 - \rho \beta + \frac12 \sqrt{1 + 4 \rho^2 \kappa} \right) \| \mathbf{W} \|_2
\end{align*}
for $\beta \equiv \epsilon \lambda_0 + c_0 + h_*^{2+p} \gamma \lambda_*$. 

Suppose $4\kappa > \beta^2$.  
Define $R_2 = \frac{2 \beta}{4 \kappa - \beta^2} > 0$.  
Observe that 
\begin{align*}
	\rho < \frac{2 \beta}{4 \kappa - \beta^2}
	& \implies \rho \left( 4 \kappa - \beta^2 \right) < 2 \beta \\ 
	& \implies \rho^2 \left( 4 \kappa - \beta^2 \right) < 2 \rho \beta \\ 
	& \implies 4 \rho^2 \kappa < 2 \rho \beta + \rho^2 \beta^2 \\ 
	& \implies 1 + 4 \rho^2 \kappa < 1 + 2 \rho \beta + \rho^2 \beta^2 = \left( 1 + \rho \beta \right)^2 . 
\end{align*}
Thus, for all $0 < \rho < R_2$, there holds 
\begin{align*}
	\frac12 - \rho \beta + \frac12 \sqrt{1 + 4 \rho^2 \kappa} 
	& < \frac12 - \rho \beta + \frac12 \left( 1 + \rho \beta \right) 
	= 1 - \frac12 \rho \beta . 
\end{align*}
Now suppose $4 \kappa \leq \beta^2$. 
Then, there holds 
\[
	\frac12 - \rho \beta + \frac12 \sqrt{ 1 + 4 \rho^2 \kappa} 
	\leq \frac12 - \rho \beta + \frac12 \sqrt{ 1 + \rho^2 \beta^2 }
	\leq \frac12 - \rho \beta + \frac12 \left( 1 + \rho \beta \right) 
	= 1 - \frac12 \rho \beta . 
\]
Combining both cases, we have 
\[
	\| \widehat{\mathbf{W}} \|_2 \leq (1 - \rho \beta / 2) \| \mathbf{W} \|_2 
\]
for all $\rho > 0$ with $\rho < \min \{ R_1 , R_2 \}$ and $\beta = \epsilon \lambda_0 + c_0 + h_*^{2+p} \gamma \lambda_*$.
The proof is complete.  \hfill
\end{proof}

%%%

By the contraction mapping theorem, it follows that \eqref{M_rho} has a unique fixed point.  
Thus, we have the following admissibility result.  
Note that admissibility holds for any choice $\gamma > 0$ even for the degenerate case when $\epsilon = c_0 = 0$.  
This is in contrast to the modified scheme and analysis technique in \cite{HJfiltered} that required $\gamma$ be sufficiently large.  
Lemma~\ref{contraction_uh} also yields the following $\ell^2$-stability result.  

%%%

\begin{theorem}
The finite difference scheme \eqref{RCD_fd} with \eqref{bc2a} has a unique solution for $\mathbf{b}$ constant-valued 
in \eqref{RCD}
whenever $\epsilon > 0$ or $c_0 > 0$ or $\gamma > 0$.  
\end{theorem}

%%%

\begin{theorem} \label{L2stable}
Suppose $\epsilon > 0$ or $c_0 = \inf_{\Omega} c > 0$ and $\mathbf{b}$ is constant-valued in \eqref{RCD}.  
Let $U$ be the solution to the FD method \eqref{RCD_fd} with \eqref{bc2a}. 
Then 
\[
	\left( \prod_{i=1,2,\ldots,d} h_i^{1/2} \right) \| U \|_{\ell^2(\cT_{\bh} \cap \Omega)} \leq C 
\]
for some constant $C$ that depends on $\Omega$, $f$, and $g$. 
\end{theorem}

\begin{proof}
Define the function $v \in H^2(\Omega) \cap C^0(\overline{\Omega})$ to be the solution to
\begin{subequations} \label{aux_v_pde}
\begin{align}
       -\Delta v &= 0  \qquad\text{in } \Omega, \\
    v &= g \qquad\text{on } \partial\Omega, 
\end{align}
\end{subequations}
and define $V: \cT_{\bh}'\to \mathbb{R} $ by $V_\alpha=v(\bx_\alpha)$ for all $\bx_\alpha \in \cT_{\bh}\cap \overline{\Omega}$ 
with ghost points introduced so that the auxiliary boundary condition \eqref{bc2a} holds for $V$. 
Let $U$ be the solution to the FD method \eqref{RCD_fd} with \eqref{bc2a}. 
Then, $U_\alpha - V_\alpha$ is the solution to 
\begin{align*}
    \widehat{L}_{\bh}^\epsilon [U_\alpha -V_\alpha] & = f(\bx_\alpha) - \widehat{L}_{\mathbf{h}}^\epsilon [V_\alpha]  
    		\quad && \text{if } \bx_\alpha \in \cT_{\bh} \cap  \Omega , \\
    U_\alpha -V_\alpha & = 0  \quad && \text{if } \bx_\alpha \in \cT_{\bh} \cap \partial \Omega , \\
    \delta_{x_i, h_i}^2 (U_\alpha - V_\alpha ) & = 0 \quad && \text{if } \bx_\alpha \in \mathcal{S}_{h_i} 
\end{align*}
by the linearity of $\widehat{L}_{\bh}^\epsilon$ and the boundary data.  

Since $U$ is a solution to \eqref{RCD_fd} and \eqref{bc2a}, by Lemma~\ref{contraction_uh} there holds 
\begin{align*}
	- \rho \| \widehat{L}_{\bh}^\epsilon V - f \|_{\ell^2(\cT_{\bh} \cap \Omega)} + \| U - V \|_{\ell^2(\cT_{\bh} \cap \Omega)} 
	& \leq \| U - V + \rho \widehat{L}_{\bh}^\epsilon V - \rho f \|_{\ell^2(\cT_{\bh} \cap \Omega)} \\ 
	& = \| \cM_\rho U - \cM_\rho V \|_{\ell^2(\cT_{\bh} \cap \Omega)} \\ 
	& \leq \left( 1 - \rho \frac{\epsilon \lambda_0 + c_0}{2} \right) \| U - V \|_{\ell^2(\cT_{\bh} \cap \Omega)} 
\end{align*}
for all $\rho > 0$ sufficiently small.  
Thus, 
\begin{equation} \label{U_Verror}
	\| U - V \|_{\ell^2(\cT_{\bh} \cap \Omega)} \leq \frac{2}{\epsilon \lambda_0 + c_0} \| \widehat{L}_{\bh}^\epsilon V - f \|_{\ell^2(\cT_{\bh} \cap \Omega)} , 
\end{equation}
and it follows that 
\[
	\| U \|_{\ell^2(\cT_{\bh} \cap \Omega)} 
	\leq \| V \|_{\ell^2(\cT_{\bh} \cap \Omega)} + \frac{2}{\epsilon \lambda_0 + c_0} \| \widehat{L}_{\bh}^\epsilon V - f \|_{\ell^2(\cT_{\bh} \cap \Omega)} . 
\] 
Note that 
\begin{align*}
	 \| \widehat{L}_{\bh}^\epsilon V \|_{\ell^2(\cT_{\bh} \cap \Omega)} 
	& \leq \left(\|\widehat{L}_{\bh}^\epsilon [v]- L^\epsilon[v] \|_{\ell^2(\cT_{\bh} \cap \Omega)} \right)
		+ \| L^\epsilon [v] \|_{\ell^2(\cT_{\bh} \cap \Omega)} , 
\end{align*}	
where $\|\widehat{L}_{\bh}^\epsilon [v]- L^\epsilon[v] \|_{\ell^2(\cT_{\bh} \cap \Omega)}$ is associated with 
the global truncation error 
for $v$ and $\| L^\epsilon[v] \|_{\ell^2(\cT_{\bh} \cap \Omega)}$ depends on the smoothness of $v$ 
which is determined by the smoothness of $g$. 
Hence, there exists a constant $C$ independent of $\bh$ such that 
\[
	\left( \prod_{i=1,2,\ldots,d} h_i^{1/2} \right) \| U \|_{\ell^2(\cT_{\bh} \cap \Omega)} \leq C .
\]
The proof is complete. 
\hfill 
\end{proof} 

For non-degenerate problems with $\epsilon > 0$ or $c_0 > 0$, we can immediately derive rates of convergence in the 
weighted $\ell^2$-norm.  

\begin{theorem}\label{conv_rates}
Suppose $u \in C^5(\overline{\Omega})$ is the solution to \eqref{RCD} with $\mathbf{b}$ constant-valued 
and $\epsilon > 0$ or $c_0 = \inf_{\Omega} c > 0$.  
Let $U$ be the solution to the FD method \eqref{RCD_fd} with \eqref{bc2a}.  
Then, 
\[
	\left( \prod_{i=1,2,\ldots,d} h_i^{1/2} \right) \| U - u \|_{\ell^2(\cT_{\bh} \cap \Omega)} \leq C h^{p+\frac12} 
\]
for some constant $C$ independent of $\bh$.  
\end{theorem}

\begin{proof}
The result follows by choosing $v$ to be the solution to \eqref{RCD} instead of \eqref{aux_v_pde} 
in the proof of Theorem~\ref{L2stable}.  
Then, by \eqref{U_Verror}, there holds 
\begin{align*}
	\| U_\alpha - u(\bx_\alpha) \|_{\ell^2(\cT_{\bh} \cap \Omega)}
	& \leq \frac{2}{ \epsilon \lambda_0 + c_0} \| \widehat{L}_{\bh}^\epsilon [u(\bx_\alpha)] - f(\bx_\alpha) \|_{\ell^2(\cT_{\bh} \cap \Omega)} \\ 
	& \leq \frac{2}{ \epsilon \lambda_0 + c_0} \| \widehat{L}_{\bh}^\epsilon [u(\bx_\alpha)] - L^\epsilon[ u(\bx_\alpha) ] \|_{\ell^2(\cT_{\bh} \cap \Omega)} \\ 
	& \qquad 
		+ \frac{2}{ \epsilon \lambda_0 + c_0} \| L^\epsilon[ u(\bx_\alpha)] - f(\bx_\alpha) \|_{\ell^2(\cT_{\bh} \cap \Omega)} \\ 
	& = \frac{2}{ \epsilon \lambda_0 + c_0} \| \widehat{L}_{\bh}^\epsilon [u(\bx_\alpha)] - L^\epsilon[ u(\bx_\alpha) ] \|_{\ell^2(\cT_{\bh} \cap \Omega)} .  
\end{align*}
The result follows by Theorem~\ref{GTE_theorem} since $\frac{2}{ \epsilon \lambda_0 + c_0}$ 
is uniformly bounded for $\epsilon > 0$ or $c_0 > 0$.  
The proof is complete. 
\hfill 
\end{proof}

\begin{remark}
The above stability and convergence analysis holds for the degenerate case $\epsilon = 0$ as long as $c_0 > 0$.  
Thus, we can scale the numerical viscosity by $h_*^2$ or even discard the numerical viscosity and maintain convergence.  
The numerical moment yielded a truncation error with order larger than one whenever choosing $p > \frac12$, 
and the analysis allows for the use of the central difference operator for the gradient exploiting the fact that it is antisymmetric.  
In practice, we know that central difference can yield inaccurate approximations on coarse meshes when $\epsilon$ is small.  
We will see in the numerical tests that adding a numerical moment improves the accuracy on course meshes while preserving 
higher-order convergence rates as $h \to 0^+$ when the underlying PDE solution is sufficiently smooth.  
\end{remark}

%%%%%%%%%%%%%%%%%%%%%%%%%%%%%%%%%%%%%%%%%%%%%%%%%%
%%%%%%%%%%%%%%%%%%%%%%%%%%%%%%%%%%%%%%%%%%%%%%%%%%
%%%%%%%%%%%%%%%%%%%%%%%%%%%%%%%%%%%%%%%%%%%%%%%%%%
%%%%%%%%%%%%%%%%%%%%%%%%%%%%%%%%%%%%%%%%%%%%%%%%%%

\section{Numerical Experiments} \label{numerics_sec}

In this section we test the convergence rates of the proposed method for approximating 
smooth solutions of \eqref{RCD}, \eqref{HJ}, and \eqref{HJeps} 
as well as the ability for the method to approximate lower regularity viscosity solutions.  
We will consider several one-dimensional and two-dimensional problems.  
Most of the tests will correspond to the nonlinear HJ equation \eqref{HJ} representing the degenerate 
case when $\epsilon = 0$.   
We will benchmark the convergence rates for both choices of the auxiliary boundary condition \eqref{bc2a} and \eqref{bc2b} 
typically choosing $\sigma = 1$ and $p=0$ or $p=1$.  
The choice $p=0$ represents the lowest-order accuracy for the proposed methods.  
We will see that the choice $p=0$ and the less accurate auxiliary boundary condition \eqref{bc2a} 
will almost always yield a more accurate approximation than the Lax-Friedrich's method 
despite the presence of a high-order boundary layer that offsets the higher-order accuracy 
associated with the central difference approximation of the gradient.  
When the underlying viscosity solution has lower regularity or when the Dirichlet boundary 
condition has a mismatch at the outflow boundary we will see that the proposed methods 
are competitive with the Lax-Friedrich's method.  
Since the outflow boundary is unknown for HJ equations in general (with the outflow boundary 
depending on the unknown solution $u$), we choose to benchmark the performance against 
the Lax-Friedrich's method even in the linear case when the advection field is known 
and upwinding methods are readily available.  

All tests will feature uniform meshes with $h_x = h_y$ for the two-dimensional experiments.  
Errors will be measured in both the weighted $\ell^2$-norm and the $\ell^\infty$-norm.   
Exploiting the uniform meshes and to allow a simple comparison with the Lax-Friedrich's method, 
we rewrite the FD scheme \eqref{HJeps_fd} as 
\begin{align} \label{scheme_general}
	\widehat{H}_{\mathbf{h}}^\epsilon [U_\alpha] 
	& \equiv - (\epsilon + \epsilon_{\bh}) \Delta_{\mathbf{h}} U_\alpha + H \left( \nabla_{\mathbf{h}} U_\alpha , U_\alpha , \bx_\alpha \right) 
		+ \gamma_{\bh} \left( \Delta_{2\mathbf{h}} - \Delta_{\mathbf{h}} \right) U_\alpha 
\end{align}
for all $\bx_\alpha \in \cT_{\bh} \cap \Omega$ 
for constants $\epsilon_{\bh} = \sigma h^r > 0$ and $\gamma_{\bh} = \gamma h^p \geq 0$ 
with $r = 1$ and $\gamma = 0$ corresponding to the Lax-Friedrich's method 
and $r = 2$ and $\gamma > 0$ corresponding to the proposed non-monotone FD method 
with a numerical moment stabilization term.  
The tests with $H$ linear correspond to the scheme \eqref{RCD_fd} for approximating \eqref{RCD}.  

The tests are all performed in {\it Matlab} and use {\it fsolve} to solve 
\eqref{scheme_general} when $H$ is nonlinear.   
To assist the nonlinear solver when $\epsilon = 0$, 
we form a sequence of $\epsilon_{\bh}$ values to generate initial guesses for {\it fsolve}. 
First we choose a constant $C$ and power $q \in \{1,2 \}$ 
so that $\epsilon = C h^1$ corresponds to the Lax-Friedrich's method.  
We then choose powers $q_k$ for $k=1,2,\ldots,N$  for 
$q = q_N > q_{N-1} > \cdots > q_1 = 0$ 
and solve \eqref{scheme_general} with an appropriate auxiliary boundary condition if $\gamma > 0$ 
using $\epsilon_{\bh} =Ch^{q_k}$ for $k=1,2,\ldots, N$.  
For $k > 1$, the initial guess for {\it fsolve} is the solution for $\epsilon_{\bh} = C h^{q_{k-1}}$.  
For $k=1$, the initial guess for {\it fsolve} is the secant line connecting the boundary data if $d=1$ 
or the zero function if $d=2$.  
We typically choose $N \leq 3$.  
The numerical tests significantly expand upon tests found in \cite{HJfiltered}.

%%%%%%%%%%%%%%%%%%%%%%%
%%%%%%%%%%%%%%%%%%%%%%%

\subsection{One-Dimensional Tests}

We first consider a series of numerical experiments in one dimension to test the accuracy of the proposed scheme. 
The first three tests are linear with the remaining test nonlinear.  
The linear problems will focus on \eqref{RCD} as the problem transitions towards convection-dominated.  
The third linear problem will focus on an example where a discontinuity at the right boundary develops 
as $\epsilon \to 0$.  
The nonlinear problem will focus on a lower-regularity solution with a corner.  
We will see that the numerical moment stabilization helps the methods perform well on coarse meshes 
when compared to the central difference method 
and yields more accurate approximations than monotone methods.  
In the tables and figures, $p=0$ and auxiliary boundary condition \eqref{bc2a} corresponds to the label ``Moment BC1", 
$p=1$ and auxiliary boundary condition \eqref{bc2a} corresponds to the label ``h*Moment BC1", 
$p=0$ and auxiliary boundary condition \eqref{bc2b} corresponds to the label ``Moment BC2", 
and $p=1$ and auxiliary boundary condition \eqref{bc2b} corresponds to the label ``h*Moment BC2".

%%%%%%%%%%%%%%%%%%%%%%%
%%%%%%%%%%%%%%%%%%%%%%%

\subsubsection{Example 1:  Problem \eqref{RCD} with a smooth solution and $b$ positive} \label{1dtest1_sec}
Consider the linear problem \eqref{RCD} with 
\[
	\epsilon = \text{1.0e-11}, \qquad b(x) = 16 \sqrt{x} + \sin(\pi x) + 2, \qquad c(x) = 0, \qquad \Omega = (0,1)
\]
with $f$ and $g$ chosen such that $u(x) = (x^2+1)\cos(\pi x) e^x$.  
The results can be found in Table~\ref{1D1_rates}.  
Notice that the central difference method is superconvergent but starts with extremely high errors for coarse meshes.  
The methods stabilized with a numerical moment have much smaller errors on coarse meshes and achieve second 
order convergence when $p=1$.  
All of the stabilized approximations are more accurate than the monotone approximations.  
Similar results hold for the $\ell^\infty$ errors. 
The small value for $\epsilon$ and the choice of $b$ ensure this problem represents 
the convection-dominated regime in a rescaled problem (\cite{Birkhoff_Gartland_Lynch90}).  

\begin{table}[htb] 
{\small
\begin{center}
\begin{tabular}{| c | c | c | c | c | c | c |}
		\hline
	 & \multicolumn{2}{|c|}{Upwind} & \multicolumn{2}{|c|}{Lax-Friedrichs} & \multicolumn{2}{|c|}{Central} \\ 
		\hline
	 $h$ & $\ell^2$ Error & Order & $\ell^2$ Error & Order & $\ell^2$ Error & Order \\ 
		\hline
	1.67e-01 & 6.91e-01 & & 9.88e-01 & & 2.06e+09 &  \\ 
		\hline
	8.33e-02 & 3.58e-01 & 0.95 & 5.36e-01 & 0.88 & 1.31e+08 & 3.98 \\ 
		\hline
	4.55e-02 & 1.99e-01 & 0.97 & 3.02e-01 & 0.95 & 1.15e+07 & 4.01 \\ 
		\hline
	1.92e-02 & 8.49e-02 & 0.99 & 1.31e-01 & 0.97 & 3.65e+05 & 4.01 \\ 
		\hline
	9.80e-03 & 4.35e-02 & 0.99 & 6.75e-02 & 0.98 & 2.45e+04 & 4.01 \\ 
		\hline
	3.31e-03 & 1.47e-02 & 1.00 & 2.31e-02 & 0.99 & 3.18e+02 & 4.00 \\ 
		\hline
	9.98e-04 & 4.44e-03 & 1.00 & 6.99e-03 & 0.99 & 2.62e+00 & 4.00 \\ 
		\hline
	2.00e-04 & 8.90e-04 & 1.00 & 1.41e-03 & 1.00 & 4.22e-03 & 4.00 \\ 
		\hline
	5.00e-05 & 2.22e-04 & 1.00 & 3.52e-04 & 1.00 & 1.65e-05 & 4.00 \\ 
		\hline
	2.00e-05 & 8.90e-05 & 1.00 & 1.41e-04 & 1.00 & 4.22e-07 & 4.00 \\ 
		\hline
\end{tabular}
\end{center}
}

\medskip 

{\small
\begin{center}
\begin{tabular}{| c | c | c | c | c | c | c | c | c |}
		\hline
	 & \multicolumn{2}{|c|}{Moment BC1} & \multicolumn{2}{|c|}{h*Moment BC1} & \multicolumn{2}{|c|}{Moment BC2} & \multicolumn{2}{|c|}{h*Moment BC2} \\ 
		\hline
	 $h$ & $\ell^2$ Error & Order & $\ell^2$ Error & Order & $\ell^2$ Error & Order & $\ell^2$ Error & Order \\ 
		\hline
	1.67e-01 & 2.67e-01 & & 2.52e-01 & & 2.49e-01 & & 2.54e-01 &  \\ 
		\hline
	8.33e-02 & 7.24e-02 & 1.88 & 6.13e-02 & 2.04 & 6.17e-02 & 2.01 & 6.16e-02 & 2.04 \\ 
		\hline
	4.55e-02 & 2.56e-02 & 1.72 & 1.75e-02 & 2.07 & 1.80e-02 & 2.03 & 1.74e-02 & 2.09 \\ 
		\hline
	1.92e-02 & 6.77e-03 & 1.55 & 2.98e-03 & 2.06 & 3.14e-03 & 2.03 & 2.87e-03 & 2.09 \\ 
		\hline
	9.80e-03 & 2.60e-03 & 1.42 & 7.57e-04 & 2.03 & 7.95e-04 & 2.04 & 7.10e-04 & 2.07 \\ 
		\hline
	3.31e-03 & 6.11e-04 & 1.34 & 8.53e-05 & 2.01 & 8.74e-05 & 2.03 & 7.76e-05 & 2.04 \\ 
		\hline
	9.98e-04 & 1.30e-04 & 1.29 & 7.77e-06 & 2.00 & 7.67e-06 & 2.03 & 6.93e-06 & 2.01 \\ 
		\hline
	2.00e-04 & 1.65e-05 & 1.28 & 3.14e-07 & 1.99 & 2.97e-07 & 2.02 & 2.76e-07 & 2.00 \\ 
		\hline
	5.00e-05 & 2.76e-06 & 1.29 & 1.97e-08 & 2.00 & 1.81e-08 & 2.02 & 1.72e-08 & 2.00 \\ 
		\hline
	2.00e-05 & 8.41e-07 & 1.30 & 3.17e-09 & 2.00 & 2.91e-09 & 1.99 & 2.76e-09 & 2.00 \\ 
		\hline
\end{tabular}
\end{center}
}
\caption{
Rates of convergence for various approximation methods for Example 1 in one dimension 
using $\epsilon_h = 9 h^q$ and $\gamma = h^p$.}
\label{1D1_rates}
\end{table}

%%%%%%%%%%%%%%%%%%%%%%%
%%%%%%%%%%%%%%%%%%%%%%%

\subsubsection{Example 2:  Problem \eqref{RCD} with a smooth solution and $b$ sign-changing} \label{1dtest2_sec}
Consider the linear problem \eqref{RCD} with 
\[
	b(x) = 16 x - 8 + \sin(\pi x), \qquad c(x) = \max\{ x-0.5, 0\} + 0.00001, \qquad \Omega = (0,1)
\]
with $f$ and $g$ chosen such that $u(x) = (x^2+1)\cos(\pi x) e^x$ for each given value of $\epsilon > 0$.  
The reaction term $c \geq \text{0.1e-4}$ is to ensure the problem satisfies a comparison principle in the limit $\epsilon = 0$.  
The problem is convection dominated throughout at least half of the domain and has a change of direction due to the sign change 
in $b$.  
The results can be found in Table~\ref{1D2a_rates} for $\epsilon = \text{1.0e-1}$, 
Table~\ref{1D2b_rates} for $\epsilon = \text{1.0e-3}$, Table~\ref{1D2c_rates} for $\epsilon = \text{1.0e-5}$, 
and Table~\ref{1D2c_rates} for $\epsilon = \text{1.0e-7}$.
Overall the central difference method performs the best with several methods having large errors on the coarse meshes 
before exhibiting behavior more consistent with having either first or second order convergence.  
The Lax-Friedrich's method and numerical moment stabilizer have less consistent rates as the stabilization terms 
are resolved.  
When comparing the actual errors, the choice $p=0$ appears more optimal for approximating on a coarse mesh 
while the choice $p=1$ appears more optimal on a finer mesh.  
Note that for the degenerate case $c=0$ and small $\epsilon$ 
all methods except those with a numerical moment and $p=0$ reported an ill-conditioned matrix.  
When letting $c=1$, the monotone methods exhibited first order convergence and all other methods exhibited second order convergence 
independent of the choice for $\epsilon$.   
The errors are all measured in $\ell^\infty$ which strongly emphasizes the potential boundary layer 
for numerical moment stabilization of the central difference method.  

\begin{table}[htb] 

{\small
\begin{center}
\begin{tabular}{| c | c | c | c | c | c | c |}
		\hline
	 & \multicolumn{2}{|c|}{Upwind} & \multicolumn{2}{|c|}{Lax-Friedrichs} & \multicolumn{2}{|c|}{Central} \\ 
		\hline
	 $h$ & $\ell^\infty$ Error & Order & $\ell^\infty$ Error & Order & $\ell^\infty$ Error & Order \\ 
		\hline
	1.67e-01 & 1.68e+01 & & 3.30e+00 & & 1.55e-01 &  \\ 
		\hline
	8.33e-02 & 3.38e+01 & -1.01 & 5.30e+00 & -0.69 & 9.65e-01 & -2.64 \\ 
		\hline
	4.55e-02 & 3.62e+01 & -0.11 & 1.11e+01 & -1.22 & 3.99e-01 & 1.46 \\ 
		\hline
	1.92e-02 & 1.90e+01 & 0.75 & 4.92e+01 & -1.73 & 5.36e-02 & 2.33 \\ 
		\hline
	9.80e-03 & 1.01e+01 & 0.94 & 8.08e+01 & -0.74 & 1.38e-02 & 2.01 \\ 
		\hline
	3.31e-03 & 3.50e+00 & 0.98 & 4.67e+01 & 0.50 & 1.57e-03 & 2.00 \\ 
		\hline
	9.98e-04 & 1.06e+00 & 0.99 & 1.73e+01 & 0.83 & 1.43e-04 & 2.00 \\ 
		\hline
	2.00e-04 & 2.14e-01 & 1.00 & 3.76e+00 & 0.95 & 5.71e-06 & 2.00 \\ 
		\hline
	5.00e-05 & 5.36e-02 & 1.00 & 9.56e-01 & 0.99 & 3.13e-07 & 2.10 \\ 
		\hline
	2.00e-05 & 2.14e-02 & 1.00 & 3.84e-01 & 1.00 & 5.45e-07 & -0.61 \\ 
		\hline
	1.00e-05 & 1.07e-02 & 1.00 & 1.92e-01 & 1.00 & 5.34e-07 & 0.03 \\ 
		\hline
\end{tabular}
\end{center}
}

\medskip

{\small
\begin{center}
\begin{tabular}{| c | c | c | c | c | c | c | c | c |}
		\hline
	 & \multicolumn{2}{|c|}{Moment BC1} & \multicolumn{2}{|c|}{h*Moment BC1} & \multicolumn{2}{|c|}{Moment BC2} & \multicolumn{2}{|c|}{h*Moment BC2} \\ 
		\hline
	 $h$ & $\ell^\infty$ Error & Order & $\ell^\infty$ Error & Order & $\ell^\infty$ Error & Order & $\ell^\infty$ Error & Order \\ 
		\hline
	1.67e-01 & 5.29e+00 & & 8.95e+00 & & 2.38e+00 & & 6.87e+00 &  \\ 
		\hline
	8.33e-02 & 1.80e+00 & 1.55 & 4.72e+01 & -2.40 & 2.04e+00 & 0.22 & 3.70e+01 & -2.43 \\ 
		\hline
	4.55e-02 & 3.34e+00 & -1.02 & 3.78e+01 & 0.37 & 1.61e+00 & 0.39 & 3.97e+01 & -0.12 \\ 
		\hline
	1.92e-02 & 7.77e+00 & -0.98 & 6.99e+00 & 1.96 & 1.40e+01 & -2.51 & 6.99e+00 & 2.02 \\ 
		\hline
	9.80e-03 & 2.37e+00 & 1.77 & 1.85e+00 & 1.97 & 2.52e+00 & 2.54 & 1.85e+00 & 1.97 \\ 
		\hline
	3.31e-03 & 2.58e-01 & 2.04 & 2.12e-01 & 1.99 & 2.58e-01 & 2.10 & 2.12e-01 & 1.99 \\ 
		\hline
	9.98e-04 & 2.35e-02 & 2.00 & 1.93e-02 & 2.00 & 2.35e-02 & 2.00 & 1.93e-02 & 2.00 \\ 
		\hline
	2.00e-04 & 9.42e-04 & 2.00 & 7.74e-04 & 2.00 & 9.42e-04 & 2.00 & 7.74e-04 & 2.00 \\ 
		\hline
	5.00e-05 & 5.78e-05 & 2.01 & 4.84e-05 & 2.00 & 5.78e-05 & 2.01 & 4.84e-05 & 2.00 \\ 
		\hline
	2.00e-05 & 1.62e-05 & 1.39 & 8.23e-06 & 1.93 & 1.61e-05 & 1.39 & 8.23e-06 & 1.93 \\ 
		\hline
	1.00e-05 & 1.96e-05 & -0.28 & 2.41e-06 & 1.77 & 1.94e-05 & -0.27 & 2.41e-06 & 1.77 \\ 
		\hline
\end{tabular}
\end{center}
}
\caption{
Rates of convergence for various approximation methods for Example 2 in one dimension 
using $\epsilon = \text{1.0e-1}$, $\epsilon_h = 9 h^q$, and $\gamma_h = h^p$.}
\label{1D2a_rates}
\end{table}

%%%

\begin{table}[htb] 

{\small
\begin{center}
\begin{tabular}{| c | c | c | c | c | c | c |}
		\hline
	 & \multicolumn{2}{|c|}{Upwind} & \multicolumn{2}{|c|}{Lax-Friedrichs} & \multicolumn{2}{|c|}{Central} \\ 
		\hline
	 $h$ & $\ell^\infty$ Error & Order & $\ell^\infty$ Error & Order & $\ell^\infty$ Error & Order \\ 
		\hline
	1.67e-01 & 1.09e+04 & & 3.64e+00 & & 3.96e+00 &  \\ 
		\hline
	8.33e-02 & 3.99e+03 & 1.46 & 7.03e+00 & -0.95 & 5.91e-01 & 2.75 \\ 
		\hline
	4.55e-02 & 9.84e+03 & -1.49 & 2.36e+01 & -2.00 & 8.38e-01 & -0.58 \\ 
		\hline
	1.92e-02 & 3.02e+03 & 1.37 & 1.87e+02 & -2.40 & 3.23e-03 & 6.46 \\ 
		\hline
	9.80e-03 & 1.35e+03 & 1.19 & 2.18e+02 & -0.23 & 8.42e-03 & -1.42 \\ 
		\hline
	3.31e-03 & 3.97e+02 & 1.13 & 4.08e+02 & -0.58 & 2.64e-02 & -1.05 \\ 
		\hline
	9.98e-04 & 1.14e+02 & 1.04 & 3.13e+03 & -1.70 & 2.07e-03 & 2.13 \\ 
		\hline
	2.00e-04 & 2.25e+01 & 1.01 & 3.92e+03 & -0.14 & 5.20e-05 & 2.29 \\ 
		\hline
	5.00e-05 & 5.61e+00 & 1.00 & 9.81e+02 & 1.00 & 1.76e-06 & 2.44 \\ 
		\hline
	2.00e-05 & 2.24e+00 & 1.00 & 3.93e+02 & 1.00 & 1.29e-05 & -2.18 \\ 
		\hline
	1.00e-05 & 1.12e+00 & 1.00 & 1.96e+02 & 1.00 & 4.64e-05 & -1.84 \\ 
		\hline
\end{tabular}
\end{center}
}

\medskip

{\small
\begin{center}
\begin{tabular}{| c | c | c | c | c | c | c | c | c |}
		\hline
	 & \multicolumn{2}{|c|}{Moment BC1} & \multicolumn{2}{|c|}{h*Moment BC1} & \multicolumn{2}{|c|}{Moment BC2} & \multicolumn{2}{|c|}{h*Moment BC2} \\ 
		\hline
	 $h$ & $\ell^\infty$ Error & Order & $\ell^\infty$ Error & Order & $\ell^\infty$ Error & Order & $\ell^\infty$ Error & Order \\ 
		\hline
	1.67e-01 & 3.33e+00 & & 6.20e+00 & & 1.74e+00 & & 5.25e+00 &  \\ 
		\hline
	8.33e-02 & 3.62e+00 & -0.12 & 1.52e+02 & -4.62 & 6.67e+03 & -11.90 & 7.78e+01 & -3.89 \\ 
		\hline
	4.55e-02 & 1.46e+00 & 1.49 & 2.30e+03 & -4.48 & 2.82e+00 & 12.82 & 1.73e+04 & -8.92 \\ 
		\hline
	1.92e-02 & 3.74e+00 & -1.09 & 6.42e+02 & 1.48 & 2.06e+00 & 0.36 & 7.07e+02 & 3.72 \\ 
		\hline
	9.80e-03 & 2.08e+01 & -2.54 & 1.39e+03 & -1.14 & 6.21e+01 & -5.06 & 1.53e+03 & -1.14 \\ 
		\hline
	3.31e-03 & 1.20e+01 & 0.51 & 3.64e+02 & 1.23 & 3.09e+01 & 0.64 & 3.98e+02 & 1.24 \\ 
		\hline
	9.98e-04 & 1.37e+01 & -0.11 & 2.91e+01 & 2.11 & 4.03e+01 & -0.22 & 3.12e+01 & 2.12 \\ 
		\hline
	2.00e-04 & 1.18e+00 & 1.53 & 7.88e-01 & 2.25 & 3.14e+00 & 1.59 & 7.89e-01 & 2.29 \\ 
		\hline
	5.00e-05 & 7.09e-02 & 2.03 & 4.91e-02 & 2.00 & 1.26e-01 & 2.32 & 4.91e-02 & 2.00 \\ 
		\hline
	2.00e-05 & 6.93e-03 & 2.54 & 7.84e-03 & 2.00 & 8.72e-03 & 2.91 & 7.84e-03 & 2.00 \\ 
		\hline
	1.00e-05 & 3.14e-02 & -2.18 & 1.92e-03 & 2.03 & 3.45e-02 & -1.98 & 1.92e-03 & 2.03 \\ 
		\hline
\end{tabular}
\end{center}
}
\caption{
Rates of convergence for various approximation methods for Example 2 in one dimension 
using $\epsilon = \text{1.0e-3}$, $\epsilon_h = 9 h^q$, and $\gamma_h = h^p$.}
\label{1D2b_rates}
\end{table}

%%%

\begin{table}[htb] 

{\small
\begin{center}
\begin{tabular}{| c | c | c | c | c | c | c |}
		\hline
	 & \multicolumn{2}{|c|}{Upwind} & \multicolumn{2}{|c|}{Lax-Friedrichs} & \multicolumn{2}{|c|}{Central} \\ 
		\hline
	 $h$ & $\ell^\infty$ Error & Order & $\ell^\infty$ Error & Order & $\ell^\infty$ Error & Order \\ 
		\hline
	1.67e-01 & 1.90e+05 & & 3.64e+00 & & 8.67e+00 &  \\ 
		\hline
	8.33e-02 & 4.48e+04 & 2.09 & 7.06e+00 & -0.95 & 1.88e+00 & 2.21 \\ 
		\hline
	4.55e-02 & 1.68e+04 & 1.61 & 2.38e+01 & -2.01 & 5.90e-01 & 1.91 \\ 
		\hline
	1.92e-02 & 2.06e+03 & 2.45 & 1.88e+02 & -2.40 & 4.22e-02 & 3.07 \\ 
		\hline
	9.80e-03 & 6.73e+02 & 1.66 & 2.21e+02 & -0.24 & 5.72e-03 & 2.97 \\ 
		\hline
	3.31e-03 & 1.06e+02 & 1.70 & 4.37e+02 & -0.63 & 5.34e-04 & 2.18 \\ 
		\hline
	9.98e-04 & 1.47e+01 & 1.65 & 4.57e+03 & -1.96 & 4.02e-05 & 2.16 \\ 
		\hline
	2.00e-04 & 2.32e+00 & 1.15 & 3.92e+03 & 0.10 & 1.70e-06 & 1.97 \\ 
		\hline
	5.00e-05 & 5.66e-01 & 1.02 & 9.82e+02 & 1.00 & 9.18e-08 & 2.10 \\ 
		\hline
	2.00e-05 & 2.25e-01 & 1.01 & 3.93e+02 & 1.00 & 2.34e-07 & -1.02 \\ 
		\hline
	1.00e-05 & 1.12e-01 & 1.00 & 1.96e+02 & 1.00 & 5.80e-07 & -1.31 \\ 
		\hline
\end{tabular}
\end{center}
}

\medskip

{\small
\begin{center}
\begin{tabular}{| c | c | c | c | c | c | c | c | c |}
		\hline
	 & \multicolumn{2}{|c|}{Moment BC1} & \multicolumn{2}{|c|}{h*Moment BC1} & \multicolumn{2}{|c|}{Moment BC2} & \multicolumn{2}{|c|}{h*Moment BC2} \\ 
		\hline
	 $h$ & $\ell^\infty$ Error & Order & $\ell^\infty$ Error & Order & $\ell^\infty$ Error & Order & $\ell^\infty$ Error & Order \\ 
		\hline
	1.67e-01 & 3.32e+00 & & 6.19e+00 & & 1.74e+00 & & 5.25e+00 &  \\ 
		\hline
	8.33e-02 & 3.67e+00 & -0.15 & 1.25e+02 & -4.34 & 1.37e+02 & -6.30 & 6.90e+01 & -3.72 \\ 
		\hline
	4.55e-02 & 1.47e+00 & 1.51 & 3.56e+03 & -5.52 & 2.97e+00 & 6.32 & 5.70e+03 & -7.28 \\ 
		\hline
	1.92e-02 & 3.57e+00 & -1.03 & 8.23e+02 & 1.70 & 2.03e+00 & 0.44 & 9.09e+02 & 2.13 \\ 
		\hline
	9.80e-03 & 1.95e+01 & -2.52 & 1.90e+03 & -1.24 & 4.64e+01 & -4.65 & 2.10e+03 & -1.24 \\ 
		\hline
	3.31e-03 & 1.49e+01 & 0.25 & 3.87e+02 & 1.47 & 4.08e+01 & 0.12 & 4.28e+02 & 1.47 \\ 
		\hline
	9.98e-04 & 2.64e+01 & -0.48 & 3.52e+01 & 2.00 & 9.47e+01 & -0.70 & 3.89e+01 & 2.00 \\ 
		\hline
	2.00e-04 & 1.22e+00 & 1.91 & 1.40e+00 & 2.00 & 6.60e+00 & 1.66 & 1.55e+00 & 2.01 \\ 
		\hline
	5.00e-05 & 7.76e-02 & 1.99 & 8.51e-02 & 2.02 & 5.98e-01 & 1.73 & 9.35e-02 & 2.03 \\ 
		\hline
	2.00e-05 & 1.85e-02 & 1.56 & 1.28e-02 & 2.06 & 1.63e-01 & 1.42 & 1.40e-02 & 2.08 \\ 
		\hline
	1.00e-05 & 1.36e-02 & 0.45 & 2.93e-03 & 2.13 & 1.30e-01 & 0.32 & 3.14e-03 & 2.15 \\ 
		\hline
\end{tabular}
\end{center}
}
\caption{
Rates of convergence for various approximation methods for Example 2 in one dimension 
using $\epsilon = \text{1.0e-5}$, $\epsilon_h = 9 h^q$, and $\gamma_h = h^p$.}
\label{1D2c_rates}
\end{table}

%%%

\begin{table}[htb]

{\small
\begin{center}
\begin{tabular}{| c | c | c | c | c | c | c |}
		\hline
	 & \multicolumn{2}{|c|}{Upwind} & \multicolumn{2}{|c|}{Lax-Friedrichs} & \multicolumn{2}{|c|}{Central} \\ 
		\hline
	 $h$ & $\ell^\infty$ Error & Order & $\ell^\infty$ Error & Order & $\ell^\infty$ Error & Order \\ 
		\hline
	1.67e-01 & 2.24e+05 & & 3.64e+00 & & 8.77e+00 &  \\ 
		\hline
	8.33e-02 & 4.99e+04 & 2.17 & 7.06e+00 & -0.95 & 1.79e+00 & 2.29 \\ 
		\hline
	4.55e-02 & 1.68e+04 & 1.79 & 2.38e+01 & -2.01 & 5.71e-01 & 1.89 \\ 
		\hline
	1.92e-02 & 2.04e+03 & 2.45 & 1.88e+02 & -2.40 & 9.88e-02 & 2.04 \\ 
		\hline
	9.80e-03 & 6.63e+02 & 1.67 & 2.21e+02 & -0.24 & 2.52e-02 & 2.03 \\ 
		\hline
	3.31e-03 & 9.91e+01 & 1.75 & 4.37e+02 & -0.63 & 2.60e-03 & 2.09 \\ 
		\hline
	9.98e-04 & 9.05e+00 & 2.00 & 4.59e+03 & -1.96 & 1.82e-03 & 0.30 \\ 
		\hline
	2.00e-04 & 3.10e-01 & 2.10 & 3.92e+03 & 0.10 & 3.29e-07 & 5.36 \\ 
		\hline
	5.00e-05 & 5.91e-02 & 1.20 & 9.82e+02 & 1.00 & 4.51e-08 & 1.43 \\ 
		\hline
	2.00e-05 & 2.33e-02 & 1.02 & 3.93e+02 & 1.00 & 5.45e-09 & 2.30 \\ 
		\hline
	1.00e-05 & 1.14e-02 & 1.03 & 1.96e+02 & 1.00 & 2.09e-09 & 1.38 \\ 
		\hline
\end{tabular}
\end{center}
}

\medskip

{\small
\begin{center}
\begin{tabular}{| c | c | c | c | c | c | c | c | c |}
		\hline
	 & \multicolumn{2}{|c|}{Moment BC1} & \multicolumn{2}{|c|}{h*Moment BC1} & \multicolumn{2}{|c|}{Moment BC2} & \multicolumn{2}{|c|}{h*Moment BC2} \\ 
		\hline
	 $h$ & $\ell^\infty$ Error & Order & $\ell^\infty$ Error & Order & $\ell^\infty$ Error & Order & $\ell^\infty$ Error & Order \\ 
		\hline
	1.67e-01 & 3.32e+00 & & 6.19e+00 & & 1.74e+00 & & 5.25e+00 &  \\ 
		\hline
	8.33e-02 & 3.68e+00 & -0.15 & 1.25e+02 & -4.34 & 1.36e+02 & -6.28 & 6.89e+01 & -3.72 \\ 
		\hline
	4.55e-02 & 1.47e+00 & 1.51 & 3.58e+03 & -5.54 & 2.97e+00 & 6.30 & 5.62e+03 & -7.26 \\ 
		\hline
	1.92e-02 & 3.57e+00 & -1.03 & 8.26e+02 & 1.71 & 2.03e+00 & 0.44 & 9.12e+02 & 2.11 \\ 
		\hline
	9.80e-03 & 1.95e+01 & -2.52 & 1.91e+03 & -1.25 & 4.63e+01 & -4.64 & 2.11e+03 & -1.25 \\ 
		\hline
	3.31e-03 & 1.49e+01 & 0.25 & 3.87e+02 & 1.47 & 4.09e+01 & 0.11 & 4.28e+02 & 1.47 \\ 
		\hline
	9.98e-04 & 2.68e+01 & -0.49 & 3.53e+01 & 2.00 & 9.62e+01 & -0.71 & 3.90e+01 & 2.00 \\ 
		\hline
	2.00e-04 & 1.22e+00 & 1.92 & 1.42e+00 & 2.00 & 6.68e+00 & 1.66 & 1.57e+00 & 2.00 \\ 
		\hline
	5.00e-05 & 7.48e-02 & 2.02 & 8.86e-02 & 2.00 & 6.07e-01 & 1.73 & 9.81e-02 & 2.00 \\ 
		\hline
	2.00e-05 & 1.83e-02 & 1.53 & 1.42e-02 & 2.00 & 1.96e-01 & 1.24 & 1.57e-02 & 2.00 \\ 
		\hline
	1.00e-05 & 2.97e-02 & -0.69 & 3.54e-03 & 2.00 & 4.11e-01 & -1.07 & 3.92e-03 & 2.00 \\ 
		\hline
\end{tabular}
\end{center}
}
\caption{
Rates of convergence for various approximation methods for Example 2 in one dimension 
using $\epsilon = \text{1.0e-7}$, $\epsilon_h = 9 h^q$, and $\gamma_h = h^p$.}
\label{1D2d_rates}
\end{table}

%%%%%%%%%%%%%%%%%%%%%%%
%%%%%%%%%%%%%%%%%%%%%%%

\subsubsection{Example 3:  Problem \eqref{RCD} with a boundary layer as $\epsilon \to 0$} \label{1dtest3_sec}
\begin{figure}[htb] 
\begin{center}
\includegraphics[width=0.45\textwidth]{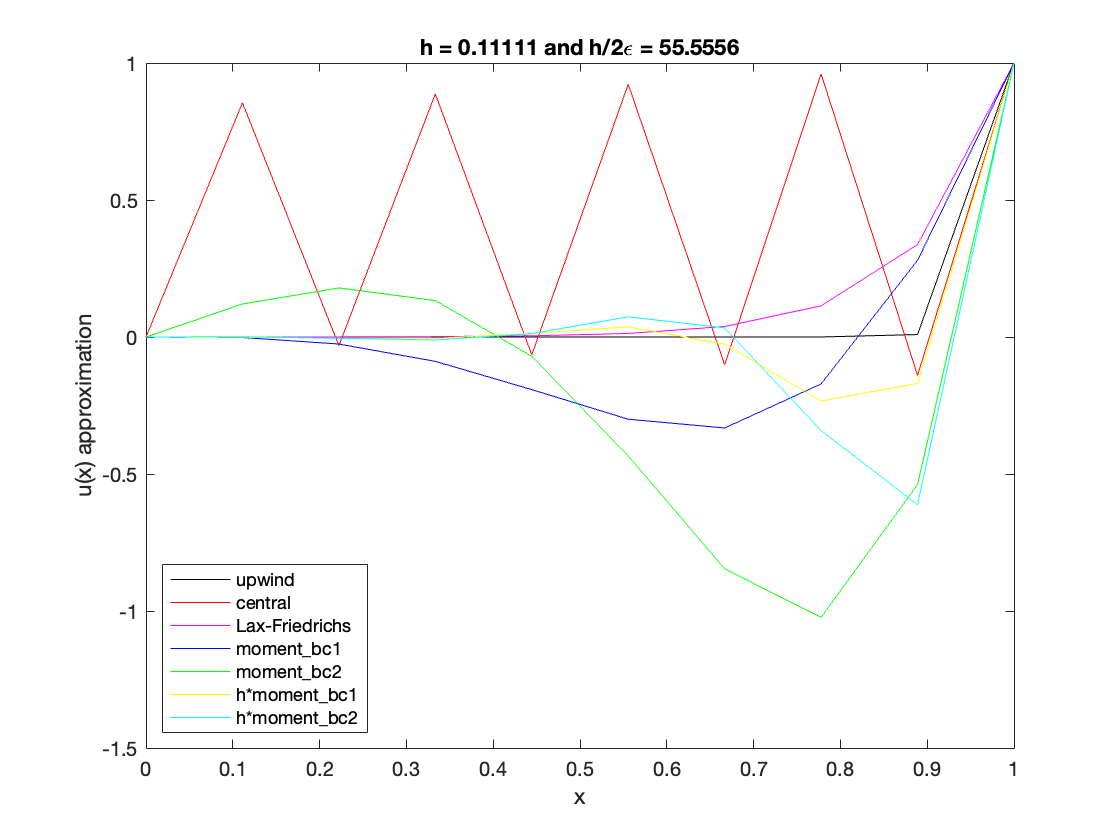} 
\qquad 
\includegraphics[width=0.45\textwidth]{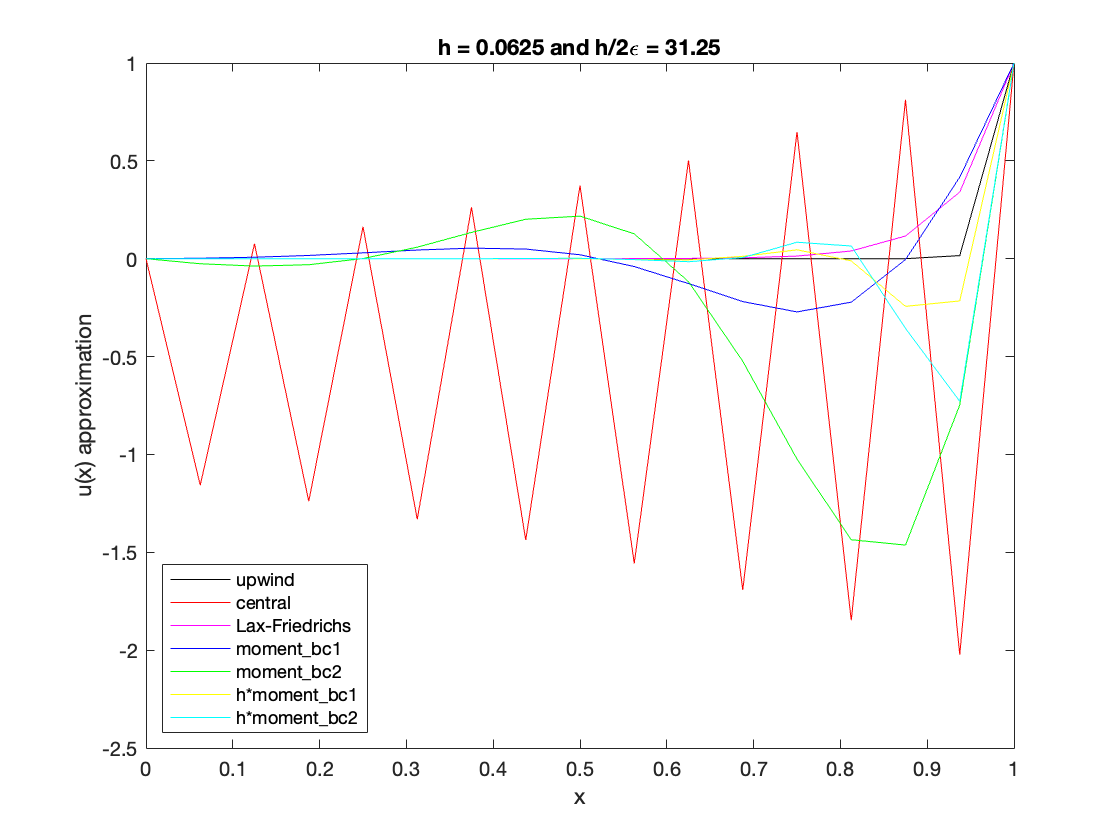} 

\includegraphics[width=0.45\textwidth]{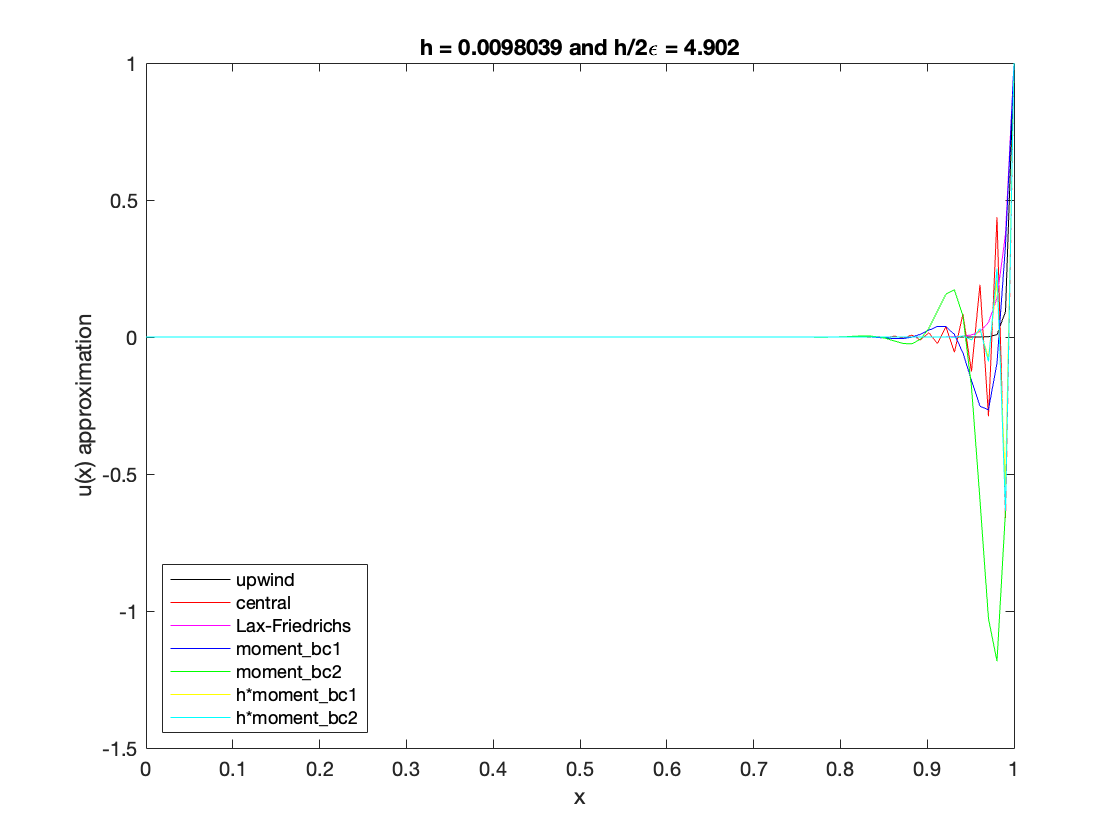} 
\qquad 
\includegraphics[width=0.45\textwidth]{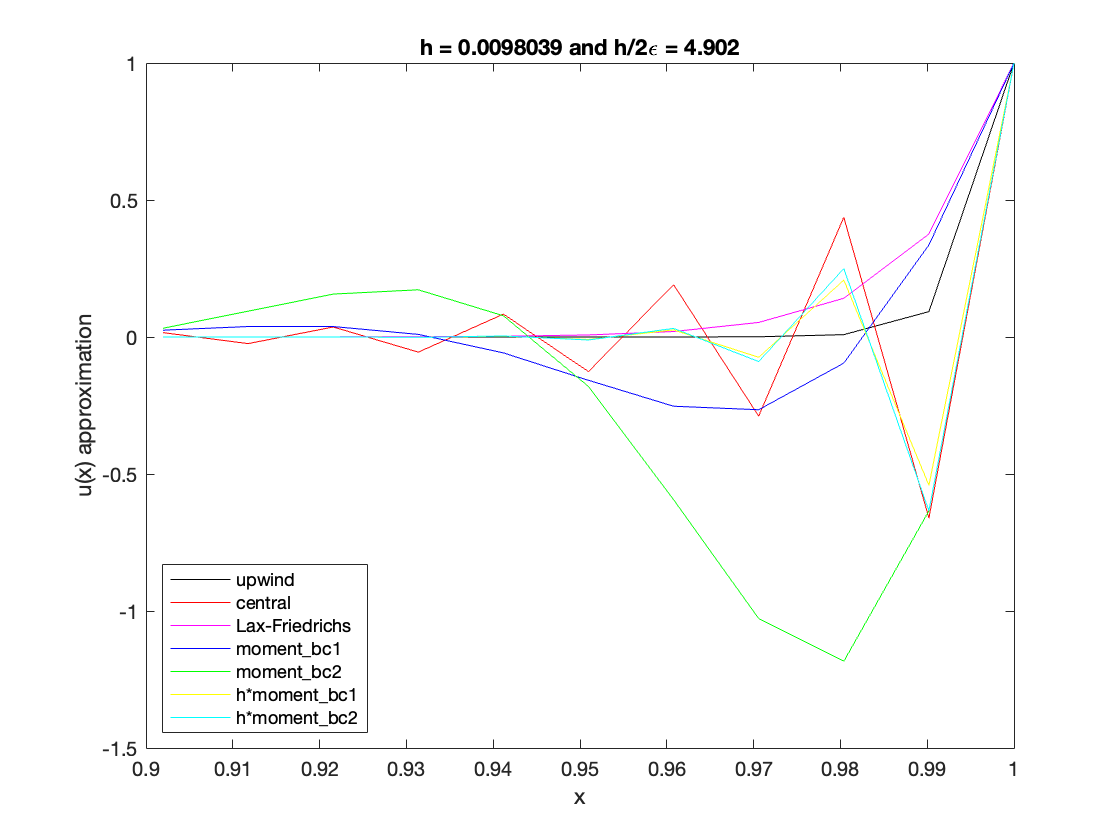} 

\includegraphics[width=0.45\textwidth]{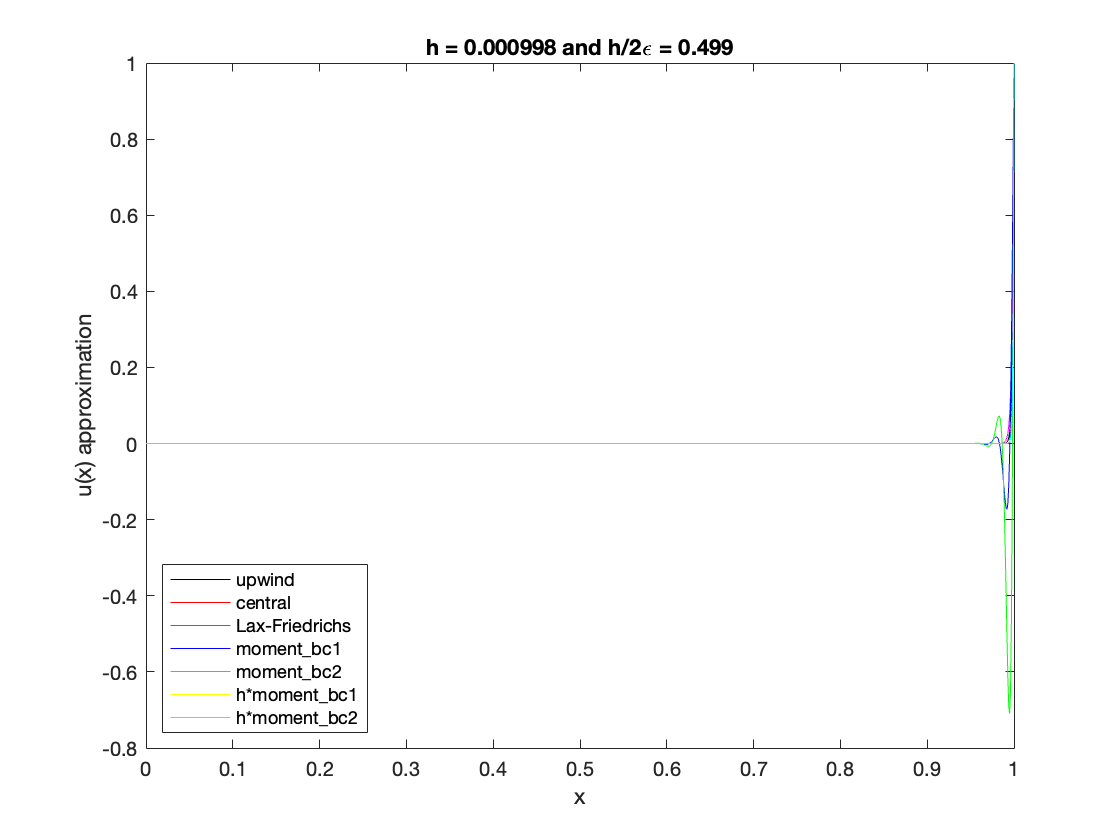} 
\qquad 
\includegraphics[width=0.45\textwidth]{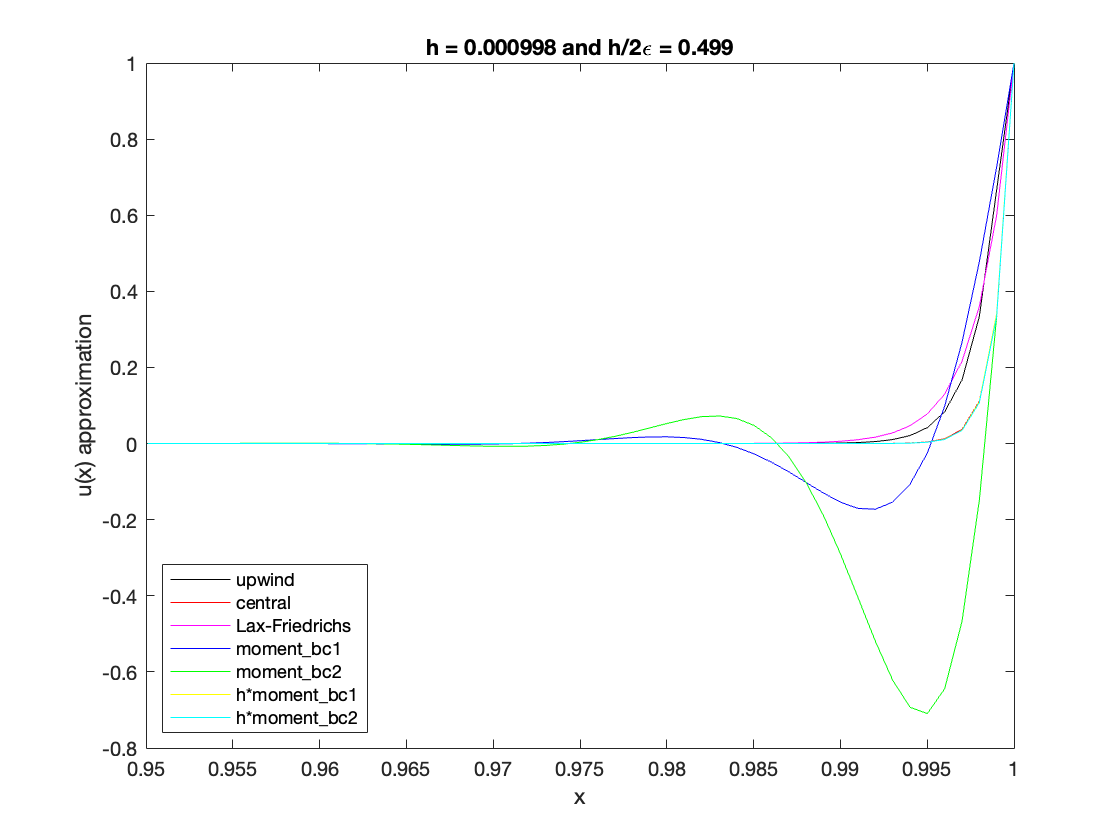}

\includegraphics[width=0.45\textwidth]{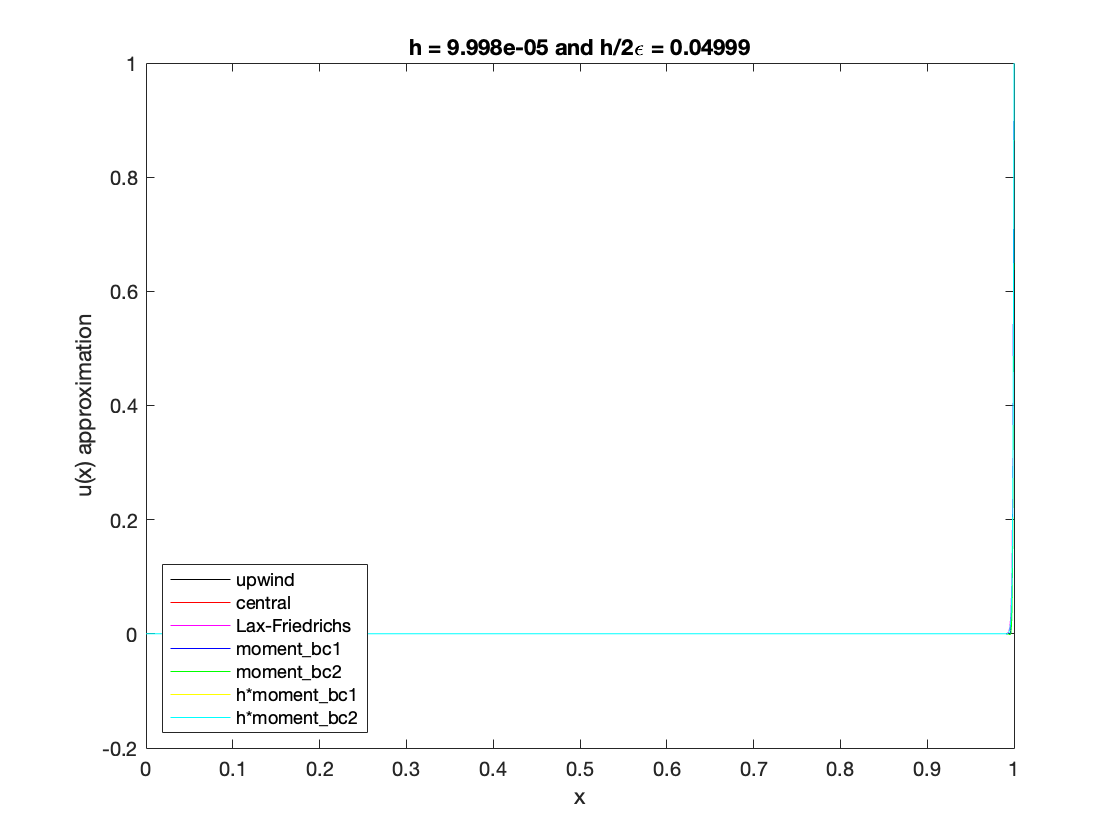} 
\qquad 
\includegraphics[width=0.45\textwidth]{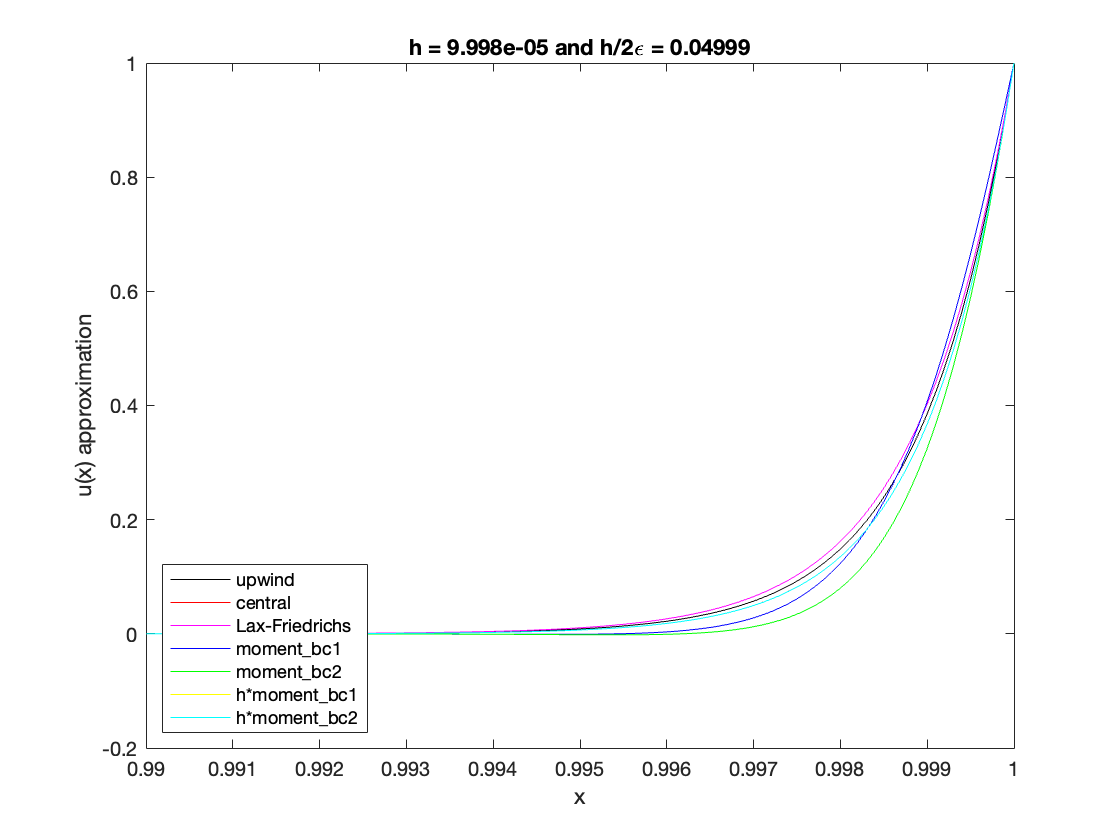} 
 
\caption{
Plots for various approximation methods and $\epsilon$ values for Example 3 in one dimension 
using $\epsilon_h = h^q$ and $\gamma_h = h^p$. 
The plots for finer meshes restrict $x$ to be near $x=1$ in the right column.  
}
\label{1D3_plots}
\end{center}
\end{figure}

Consider the linear, constant-coefficient problem \eqref{RCD} with 
\[
	b(x) = 1, \qquad c(x) = 0, \qquad f(x) = 0, \qquad \Omega = (0,1)
\]
and $g(0) = 0$, $g(1) = 1$.
The solution for $\epsilon = 0$ is $u(x) = 0$, and it satisfies the boundary data in the ``viscosity sense".  
For $0 < \epsilon \ll 1$ there is a boundary layer transition from zero to the boundary value 1.    
We can see several plots for decreasing $\epsilon$ values in Figure~\ref{1D3_plots}.  
The non-monotone methods have oscillations, but the numerical moment with $p=1$ 
suppresses the oscillations closer to $x=1$ as well as their amplitudes.  
The auxiliary boundary condition \eqref{bc2b} leads to larger errors than \eqref{bc2a}.  
We conjecture that this is due to the fact the boundary condition has a wider stencil reaching two nodes into the interior 
causing larger oscillations.  
We see that the central difference method yields extremely poor approximations until $h$ is sufficiently small ($h/2\epsilon < 1$)
to ensure the underlying method is monotone.  
Rates of convergence for $\epsilon = 0$ are recorded in Table~\ref{1D3_rates} 
where we see that the numerical moment stabilization with $p=1$ 
is competitive with the Lax-Friedrich's method despite the presence of oscillations.  
The boundary layer error for $p=0$ is much larger over a bigger region with non-convergent behavior for the auxiliary boundary 
condition \eqref{bc2b}.   
The upwind method would yield the exact answer for $\epsilon = 0$ since it would not utilize the incorrect 
boundary value at the outflow boundary and is exact for linear functions.

\begin{table}[htb]
{\small 
\begin{center}
\begin{tabular}{| c | c | c | c | c | c | c |}
		\hline
	 & \multicolumn{2}{|c|}{Lax-Friedrichs} & \multicolumn{2}{|c|}{Moment BC1} & \multicolumn{2}{|c|}{Moment BC2} \\ 
		\hline
	 $h$ & $\ell^2$ Error & Order & $\ell^2$ Error & Order & $\ell^2$ Error & Order \\ 
		\hline
	1.67e-01 & 1.44e-01 & & 2.16e-01 & & 5.21e-01 &  \\ 
		\hline
	1.11e-01 & 1.18e-01 & 0.49 & 1.99e-01 & 0.20 & 5.16e-01 & 0.02 \\ 
		\hline
	6.25e-02 & 8.84e-02 & 0.50 & 1.53e-01 & 0.46 & 6.50e-01 & -0.40 \\ 
		\hline
	4.76e-02 & 7.72e-02 & 0.50 & 1.48e-01 & 0.11 & 6.56e-01 & -0.04 \\ 
		\hline
	9.80e-03 & 3.50e-02 & 0.50 & 1.00e-01 & 0.25 & 7.01e-01 & -0.04 \\ 
		\hline
	9.98e-04 & 1.12e-02 & 0.50 & 5.17e-02 & 0.29 & 7.08e-01 & -0.00 \\ 
		\hline
	1.00e-04 & 3.54e-03 & 0.50 & 2.51e-02 & 0.31 & 7.07e-01 & 0.00 \\ 
		\hline
	1.00e-05 & 1.12e-03 & 0.50 & 1.19e-02 & 0.32 & 7.07e-01 & 0.00 \\ 
		\hline
	1.00e-06 & 3.54e-04 & 0.50 & 5.57e-03 & 0.33 & 7.07e-01 & 0.00 \\ 
		\hline
\end{tabular}
\end{center}

\medskip 

\begin{center}
\begin{tabular}{| c | c | c | c | c |}
		\hline
	 & \multicolumn{2}{|c|}{h*Moment BC1} & \multicolumn{2}{|c|}{h*Moment BC2} \\ 
		\hline
	 $h$ & $\ell^2$ Error & Order & $\ell^2$ Error & Order \\ 
		\hline
	1.67e-01 & 9.99e-02 & & 2.39e-01 &  \\ 
		\hline
	1.11e-01 & 1.00e-01 & -0.01 & 2.43e-01 & -0.04 \\ 
		\hline
	6.25e-02 & 8.61e-02 & 0.26 & 2.16e-01 & 0.20 \\ 
		\hline
	4.76e-02 & 7.85e-02 & 0.34 & 1.99e-01 & 0.31 \\ 
		\hline
	9.80e-03 & 3.97e-02 & 0.43 & 1.03e-01 & 0.41 \\ 
		\hline
	9.98e-04 & 1.30e-02 & 0.49 & 3.40e-02 & 0.49 \\ 
		\hline
	1.00e-04 & 4.12e-03 & 0.50 & 1.08e-02 & 0.50 \\ 
		\hline
	1.00e-05 & 1.30e-03 & 0.50 & 3.42e-03 & 0.50 \\ 
		\hline
	1.00e-06 & 4.12e-04 & 0.50 & 1.08e-03 & 0.50 \\ 
		\hline
\end{tabular}
\end{center}
}
\caption{
Rates of convergence for various approximation methods for Example 3 in one dimension 
using $\epsilon = 0$, $\epsilon_h = h^q$, and $\gamma_h = h^p$.}
\label{1D3_rates}
\end{table}

We now qualitatively look at the impact of the numerical moment stabilizer.  
The central difference method has high-frequency oscillations in Figure~\ref{1D3_plots} while 
the numerical moment yields oscillations with much lower frequency.  
Note that the finite difference method we are using to approximate Example 3 is given by 
\[-\sigma h^2\delta_h^2U_\alpha-\epsilon\delta_h^2U_\alpha+\delta_h U_\alpha +\gamma h^p(\delta_h^2-\delta_{2h}^2)U_\alpha =0.\]
Let $a=\sigma h^2+\epsilon$ and $b=\gamma h^p$.  
Then, the characteristic polynomial is $\rho(x)=(x-1)\eta(x)$, where $\eta(x)=bx^3 +(2h-4a-3b)x^2+ (4a+3b+2h)x-b$.  
To explore the impact of the numerical moment, we study the three roots of $\eta(x)$, namely $r_1, r_2, r_3$. 
We observe that $\eta(1)=4h\to 0$ as $h\to 0$, which means $r_1\to 1$ as $h\to 0$. 
Given that the product of roots satisfies $r_1r_2r_3=1$, there holds $r_2r_3\to 1$ as $h\to 0$.  
We focus on $r_2$ and $r_3$, and note that for the central difference method the characteristic polynomial 
has roots of $1$ and $-1$ for certain choices of parameters 
with the negative, real-valued root being responsible for the high-frequency oscillatory behavior.  
For the monotone methods, the characteristic polynomial would have two positive roots.  

Rewriting $\eta(x)$, we have 
\begin{align*}
    \eta(x)&=bx^3 +(2h-4a-3b)x^2+ (4a+3b+2h)x-b\\
    &= b\left[(x-1)^3+ x\left(-\frac{4a-2h}{b}x+\frac{4a+2h}{b} \right)\right],
\end{align*}
where the ratio $\frac{4a-2h}{b}$ plays a significant role in influencing the behavior of the roots $r_2$ and $r_3$. 
Qualitatively, 
\begin{itemize}
    \item when $\frac{4a-2h}{b} \lesssim -4$, $r_2\to 0^-$ and $r_3\to -\infty$ as $h \to 0^+$.  
    The two negative roots result in sharp osculations.
    \item when $-4 \lesssim \frac{4a-2h}{b} \lesssim 0$, $r_2$ and $r_3$ are complex-valued.  
    We observe damped wave-like oscillations.
    \item when $\frac{4a-2h}{b} \gtrsim 0$, $r_2\to 0^+$ and $r_3\to\infty$ as $h \to 0^+$. 
    The two positive roots yield exponential functions with a monotone approximation depending on the data.    
\end{itemize}
Plugging in values for $a$ and $b$,
\[
\frac{4a-2h}{b}=\frac{4(\sigma h^2+\epsilon)-2h}{\gamma h^p}.
\]
As long as $h$ is small and $\frac{\epsilon}{h^p}$ is not too small, there are two positive roots when $0\leq p<1$. 
When $p=1$, 
\[
\frac{4a-2h}{b}=\frac{4(\sigma h^2+\epsilon)-2h}{\gamma h}=\frac{2}{\gamma}\left(\frac{2\epsilon}{h}-1\right)+\mathcal{O}(h^2) , 
\]
and we have positive roots with an asymptotic mesh condition similar to that for the central difference method.  
Overall, the plots in Figure~\ref{1D3_plots} are consistent with the presence of complex roots which qualitatively 
appears to be an improvement over the central difference method for coarse meshes.  
The errors recorded in Table~\ref{1D3_rates} for the degenerate problem with $\epsilon = 0$ for coarse meshes 
support the observation that the numerical moment stabilizes the central difference method by qualitatively changing 
the approximations in a way that improves accuracy.

%%%%%%%%%%%%%%%%%%%%%%%
%%%%%%%%%%%%%%%%%%%%%%%

\subsubsection{Example 4:  Non-smooth solution with Dirichlet data} \label{1dtest4_sec}

\begin{figure}[htb] 
\begin{center}
\includegraphics[width=0.43\textwidth]{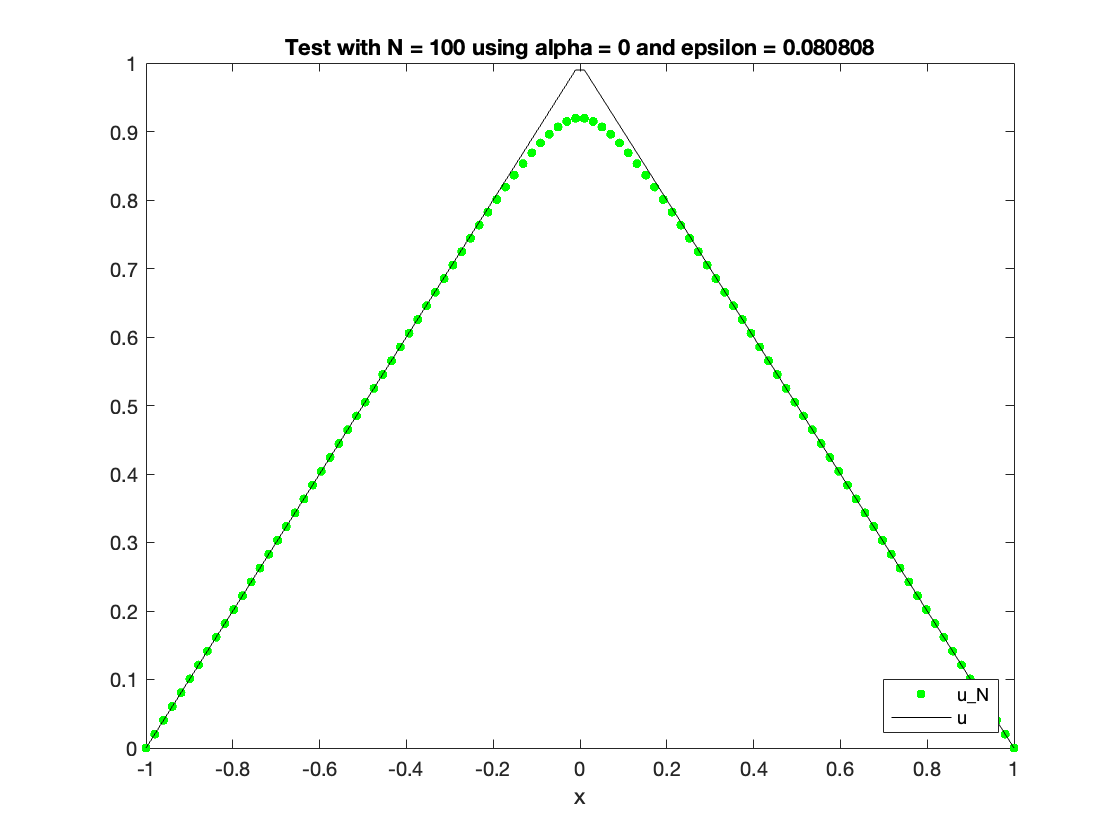} 
\qquad 
\includegraphics[width=0.43\textwidth]{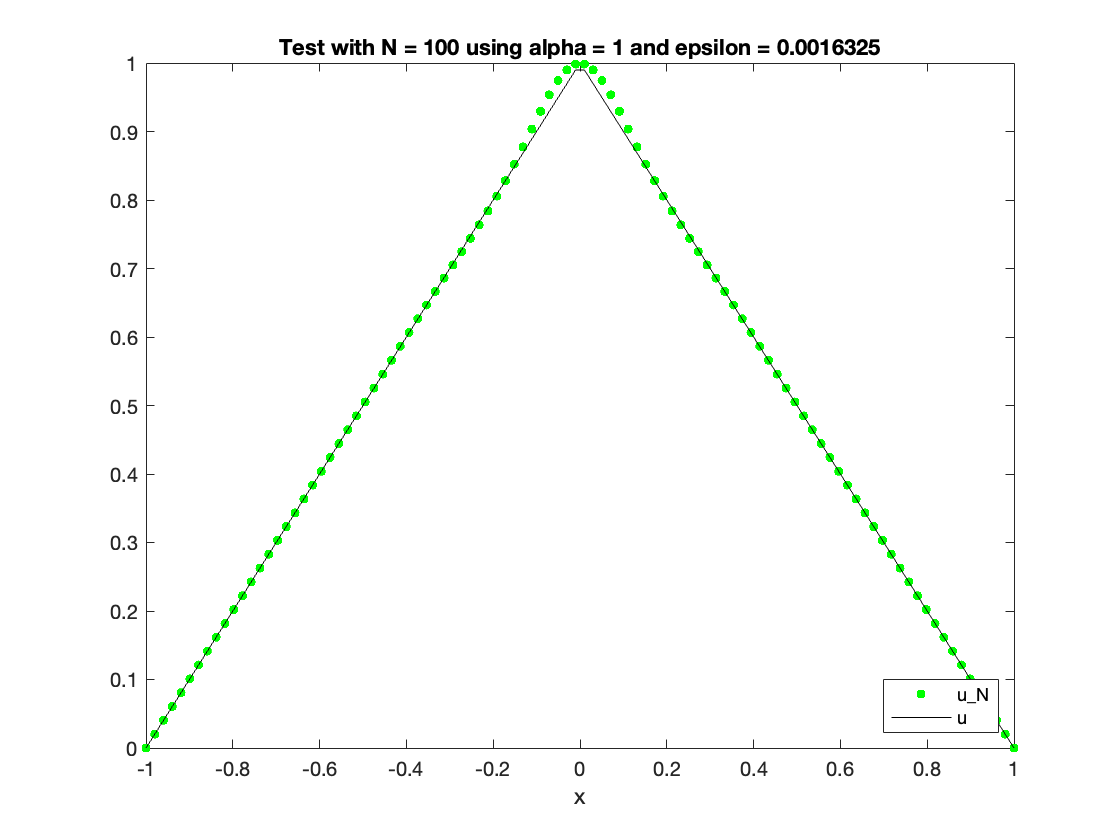} 
\caption{
Approximations for Example 4 in one dimension 
using the Lax-Friedrich's method ($\epsilon_h = 4h, \gamma_h =0$) on the left 
and the proposed non-monotone method with $\epsilon_h = 4h^2, \gamma_h = 1$ on the right.  
Both approximations correspond to $h${\em=2.02e-02}.  
}
\label{1d4_plots_bc1}
\end{center}
\end{figure}

Consider the nonlinear problem \eqref{HJ} with  
\[
	H[u] = \left| u_x \right| -1 = 0  \quad \text{in } \Omega = (-1,1)  
\]
with boundary data chosen such that the viscosity solution is $u(x) = 1 - |x|$.  
The results for the Lax-Friedrich's method and the proposed method using auxiliary boundary condition \eqref{bc2a} 
with $\gamma_h = h^p$ for $p=0,1$ can be found in Table~\ref{1D4_BC1_rates}.  
Identical results were computed for the auxiliary boundary condition \eqref{bc2b} 
with $\gamma_h = h^p$ for $p=0,1$. 
Plots of the approximations for $h$=2.02e-02 can be bound in Figure~\ref{1d4_plots_bc1}.  
We see that Lax-Friedrich's method yields a solution that is concave down over the entire domain 
like the exact viscosity solution; however,  
it has a larger error when resolving the corner of the solution.  
The proposed method switches concavity, but reduces the overall error near the corner.  
Consequently, the Lax-Friedrich's method yields a better qualitative solution 
while the proposed non-monotone method yields a slightly more accurate solution when $p=1$.  
The central difference method converged to algebraic artifacts that did not correspond to the 
viscosity solution.  
The artifacts correspond to various combinations of positive and negative signs for $u_x$ approximations 
at each node.  
Thus, the numerical moment stabilization appears to enforce uniqueness while steering the approximation 
towards the viscosity solution instead of a less-smooth artifact.  

\begin{table}[htb] 
{\small 
\begin{center}
\begin{tabular}{| c | c | c | c | c |}
		\hline 
	& \multicolumn{4}{| c |}{Lax-Friedrichs} \\ 
		\hline
	 $h$ &  $\ell^2$ Error & Rate & $\ell^\infty$ Error & Rate \\
		\hline
	2.02e-02 & 2.26e-02 &  & 7.07e-02 &  \\ 
		\hline
	6.69e-03 & 4.31e-03 & 1.50 & 2.34e-02 & 1.00  \\ 
		\hline
	3.34e-03 & 1.52e-03 & 1.50 & 1.17e-02 & 1.00  \\ 
		\hline
	2.00e-03 & 7.05e-04 & 1.50 & 7.01e-03 & 1.00  \\ 
		\hline
	1.00e-03 & 2.49e-04 & 1.50 & 3.50e-03 & 1.00 \\ 
		\hline
	5.00e-04 & 8.81e-05 & 1.50 & 1.75e-03 & 1.00 \\ 
		\hline
\end{tabular}
\end{center}

\smallskip 

\begin{center}
\begin{tabular}{| c | c | c | c | c | c | c | c | c |}
		\hline 
	& \multicolumn{4}{| c |}{Moment BC1} & \multicolumn{4}{| c |}{h*Moment BC1}\\ 
		\hline
	 $h$ &  $\ell^2$ Error & Rate & $\ell^\infty$ Error & Rate & $\ell^2$ Error & Rate & $\ell^\infty$ Error & Rate \\ 
		\hline
	2.02e-02 & 9.98e-03 &  & 2.49e-02 &  & 1.91e-03 &  & 8.47e-03 &  \\ 
		\hline
	6.69e-03 & 3.36e-03 & 0.99 & 1.21e-02 & 0.65 & 4.00e-04 & 1.42 & 3.17e-03 & 0.89 \\ 
		\hline
	3.34e-03 & 1.68e-03 & 1.00 & 7.68e-03 & 0.66 & 1.44e-04 & 1.47 & 1.62e-03 & 0.96 \\ 
		\hline
	2.00e-03 & 1.00e-03 & 1.00 & 5.46e-03 & 0.67 & 6.76e-05 & 1.48 & 9.85e-04 & 0.98 \\ 
		\hline
	1.00e-03 & 5.02e-04 & 1.00 & 3.45e-03 & 0.66 & 2.40e-05 & 1.49 & 4.96e-04 & 0.99 \\ 
		\hline
	5.00e-04 & 2.50e-04 & 1.00 & 2.17e-03 & 0.67 & 8.53e-06 & 1.50 & 2.49e-04 & 0.99 \\ 
		\hline
\end{tabular}
\end{center}

}
\caption{
Rates of convergence for Example 4 in one dimension using the Lax-Friedrich's method with $\epsilon_h = 4h, \gamma=0$ 
and the numerical moment and boundary condition \eqref{bc2a} 
with $\epsilon_h = 4h^2, \gamma_h =1$
and $\epsilon_h = 4h^2, \gamma_h =h$.  
Note that identical results were found for the numerical moment and boundary condition \eqref{bc2b} 
with $\epsilon_h = 4h^2, \gamma_h =1$
and $\epsilon_h = 4h^2, \gamma_h =h$.}
\label{1D4_BC1_rates}
\end{table}

%%%%%%%%%%%%%%%%%%%%%%%
%%%%%%%%%%%%%%%%%%%%%%%
%%%%%%%%%%%%%%%%%%%%%%%
%%%%%%%%%%%%%%%%%%%%%%%
%%%%%%%%%%%%%%%%%%%%%%%
%%%%%%%%%%%%%%%%%%%%%%%

\subsection{Two-Dimensional Tests} \label{2D_numerics_sec}

We now consider a series of numerical experiments in two dimensions to test the accuracy of the proposed scheme. 
All of the tests are performed on the domain $\Omega = (0,1)^2$ with uniform meshes where $h_x = h_y$.  
Overall, the Lax-Friedrich's method has a rate of convergence approaching 1 in the weighted $\ell^2$-norm 
and $\ell^\infty$-norm as expected, 
and it is less accurate than the proposed non-monotone FD method that adds a numerical moment stabilizer.  
However, the Lax-Friedrich's method's rates are often suboptimal with the tests appearing to be 
in the pre-asymptotic range for the $h$ values.  
The tests will all correspond to the case $\epsilon = 0$, and the Dirichlet boundary data is chosen 
to not create a boundary layer.  
The first test is linear with constant-coefficients while the 
other two tests correspond to nonlinear problems with varying degrees of regularity for the viscosity solution.  

%%%%%%%%%%%%%%%%%%%%%%%
%%%%%%%%%%%%%%%%%%%%%%%

\subsubsection{Example 1:  A Linear Benchmark} 
Consider the linear problem \eqref{RCD} with 
\begin{align*} 
	\epsilon = 0, \qquad \mathbf{b} = [1, 1]^T, \qquad c = 0, 
\end{align*}
where $f$ and $g$ chosen such that $u(x,y)=e^{xy}$. 
The test results can be found in Tables~\ref{2d1_monotone_rates}, \ref{2d1_momentBC1_rates}, and \ref{2d1_momentBC2_rates}.  
We see that the correctly scaled numerical moment allows rates of convergence near 2 in the weighted $\ell^2$-norm.  

\begin{table}[htb] 
{\small
\begin{center}
\begin{tabular}{| c | c | c | c | c | c | c | c | c |}
		\hline
	& \multicolumn{4}{|c|}{Upwind} & \multicolumn{4}{|c|}{Lax-Friedrichs} \\ 
		\hline 
	 $h$ & $\ell^2$ Error & Order & $\ell^\infty$ Error & Order & $\ell^2$ Error & Order & $\ell^\infty$ Error & Order \\
		\hline
	1.57e-01 & 1.25e-02 &  & 4.73e-02 &  & 1.94e-02 &  & 5.06e-02 &  \\
		\hline
	7.44e-02 & 6.60e-03 & 0.85 & 2.76e-02 & 0.72 & 1.15e-02 & 0.70 & 3.69e-02 & 0.42 \\ 
		\hline
	3.63e-02 & 3.37e-03 & 0.94 & 1.48e-02 & 0.87 & 6.27e-03 & 0.85 & 2.27e-02 & 0.68 \\ 
		\hline
	2.40e-02 & 2.26e-03 & 0.97 & 1.01e-02 & 0.93 & 4.30e-03 & 0.91 & 1.65e-02 & 0.77 \\ 
		\hline
	1.79e-02 & 1.70e-03 & 0.98 & 7.63e-03 & 0.95 & 3.28e-03 & 0.94 & 1.29e-02 & 0.84 \\ 
		\hline
	1.19e-02 & 1.14e-03 & 0.98 & 5.13e-03 & 0.97 & 2.22e-03 & 0.95 & 9.09e-03 & 0.85 \\ 
		\hline
	8.89e-03 & 8.52e-04 & 0.99 & 3.87e-03 & 0.98 & 1.67e-03 & 0.97 & 7.02e-03 & 0.89 \\ 
		\hline
	5.92e-03 & 5.69e-04 & 0.99 & 2.59e-03 & 0.98 & 1.12e-03 & 0.98 & 4.82e-03 & 0.92 \\ 
		\hline
\end{tabular}
\end{center}
}
\caption{
Rates of convergence for Example 1 in two dimensions using the upwinding method and the Lax-Friedrich's method 
with $\epsilon_h = h_x$.  }
\label{2d1_monotone_rates}
\end{table}

%%%

\begin{table}[htb] 
{\small
\begin{center}
\begin{tabular}{| c | c | c | c | c | c | c | c | c |}
		\hline
	& \multicolumn{4}{|c|}{Moment BC1} & \multicolumn{4}{|c|}{h*Moment BC1} \\ 
		\hline 
	 $h$ & $\ell^2$ Error & Order & $\ell^\infty$ Error & Order & $\ell^2$ Error & Order & $\ell^\infty$ Error & Order \\
		\hline
	1.57e-01 & 2.40e-02 &  & 5.05e-02 &  & 6.74e-03 &  & 2.22e-02 &  \\ 
		\hline
	7.44e-02 & 7.86e-03 & 1.50 & 2.52e-02 & 0.93 & 1.49e-03 & 2.02 & 7.96e-03 & 1.37 \\ 
		\hline
	3.63e-02 & 2.80e-03 & 1.43 & 1.26e-02 & 0.96 & 3.21e-04 & 2.13 & 2.42e-03 & 1.66 \\ 
		\hline
	2.40e-02 & 1.54e-03 & 1.44 & 8.09e-03 & 1.08 & 1.33e-04 & 2.14 & 1.15e-03 & 1.80 \\ 
		\hline
	1.79e-02 & 1.01e-03 & 1.44 & 5.93e-03 & 1.06 & 7.14e-05 & 2.12 & 6.67e-04 & 1.86 \\ 
		\hline
	1.19e-02 & 5.63e-04 & 1.44 & 3.74e-03 & 1.13 & 3.02e-05 & 2.10 & 3.06e-04 & 1.90 \\ 
		\hline
	8.89e-03 & 3.72e-04 & 1.43 & 2.67e-03 & 1.16 & 1.66e-05 & 2.08 & 1.75e-04 & 1.93 \\ 
		\hline
	5.92e-03 & 2.08e-04 & 1.43 & 1.64e-03 & 1.19 & 7.15e-06 & 2.06 & 7.92e-05 & 1.95 \\ 
		\hline
\end{tabular}
\end{center}
}
\caption{
Rates of convergence for Example 1 in two dimensions using the numerical moment and boundary condition \eqref{bc2a} 
with $\epsilon_h = h_x^2, \gamma_h = 4$ and $\epsilon_h = h_x^2, \gamma_h = h_x$.  }
\label{2d1_momentBC1_rates}
\end{table}

%%%%%

\begin{table}[htb] 
{\small
\begin{center}
\begin{tabular}{| c | c | c | c | c | c | c | c | c |}
		\hline
	& \multicolumn{4}{|c|}{Moment BC2} & \multicolumn{4}{|c|}{h*Moment BC2} \\ 
		\hline 
	 $h$ & $\ell^2$ Error & Order & $\ell^\infty$ Error & Order & $\ell^2$ Error & Order & $\ell^\infty$ Error & Order \\
		\hline
	1.57e-01 & 6.90e-03 &  & 2.01e-02 &  & 5.36e-03 &  & 1.96e-02 &  \\
		\hline
	7.44e-02 & 3.27e-03 & 1.00 & 1.37e-02 & 0.51 & 1.22e-03 & 1.98 & 8.12e-03 & 1.18 \\
		\hline
	3.63e-02 & 1.16e-03 & 1.44 & 7.65e-03 & 0.81 & 2.48e-04 & 2.21 & 2.54e-03 & 1.62 \\ 
		\hline
	2.40e-02 & 6.06e-04 & 1.58 & 5.22e-03 & 0.92 & 9.85e-05 & 2.23 & 1.21e-03 & 1.78 \\
		\hline
	1.79e-02 & 3.75e-04 & 1.64 & 3.90e-03 & 1.00 & 5.16e-05 & 2.21 & 7.08e-04 & 1.85 \\
		\hline
	1.19e-02 & 1.88e-04 & 1.69 & 2.58e-03 & 1.01 & 2.11e-05 & 2.18 & 3.26e-04 & 1.89 \\ 
		\hline
	8.89e-03 & 1.14e-04 & 1.73 & 1.89e-03 & 1.06 & 1.13e-05 & 2.15 & 1.87e-04 & 1.92 \\ 
		\hline
	5.92e-03 & 5.54e-05 & 1.76 & 1.21e-03 & 1.10 & 4.77e-06 & 2.12 & 8.46e-05 & 1.95 \\ 
		\hline
\end{tabular}
\end{center}
}
\caption{
Rates of convergence for Example 1 in two dimensions using the numerical moment and boundary condition \eqref{bc2b} 
with $\epsilon_h = h_x^2, \gamma_h = 4$ and $\epsilon_h = h_x^2, \gamma_h = h_x$.  }
\label{2d1_momentBC2_rates}
\end{table}

%%%%%%%%%%%%%%%%%%%%%%%
%%%%%%%%%%%%%%%%%%%%%%%

\subsubsection{Example 2:  A nonlinear $C^1 \backslash C^2$ operator with a smooth solution} 

Consider the nonlinear problem \eqref{HJ} with 
\begin{align*}
	H[u] \equiv \sqrt{u_x^2 + u_y^2} + u & = f \qquad \text{in } \Omega , \\ 
	u & = g \qquad  \text{on } \partial \Omega , 
\end{align*}
where $f$ and $g$ are chosen such that the exact solution is $u(x,y)=e^{xy}$.
The results can be found in Tables~\ref{2d2_monotone_rates}, \ref{2d2_momentBC1_rates}, and \ref{2d2_momentBC2_rates}.  
We see that the correctly scaled numerical moment allows rates of convergence near 2 in the weighted $\ell^2$-norm.  

\begin{table}[htb] 
{\small
\begin{center}
\begin{tabular}{| c | c | c | c | c |}
		\hline
	& \multicolumn{4}{|c|}{Lax-Friedrichs} \\ 
		\hline 
	 $h$ & $\ell^2$ Error & Order & $\ell^\infty$ Error & Order \\
		\hline
	1.57e-01 & 1.75e-02 &  & 4.47e-02 &  \\ 
		\hline
	7.44e-02 & 1.11e-02 & 0.61 & 3.27e-02 & 0.42 \\ 
		\hline
	3.63e-02 & 6.46e-03 & 0.76 & 2.16e-02 & 0.58 \\ 
		\hline
	2.40e-02 & 4.56e-03 & 0.84 & 1.63e-02 & 0.68 \\ 
		\hline
	1.79e-02 & 3.53e-03 & 0.88 & 1.30e-02 & 0.77 \\ 
		\hline
	1.19e-02 & 2.43e-03 & 0.91 & 9.37e-03 & 0.80 \\ 
		\hline
	8.89e-03 & 1.85e-03 & 0.93 & 7.30e-03 & 0.86 \\ 
		\hline
	5.92e-03 & 1.26e-03 & 0.95 & 5.09e-03 & 0.88 \\ 
		\hline
\end{tabular}
\end{center}
}
\caption{
Rates of convergence for Example 2 in two dimensions using the Lax-Friedrich's method 
with $\epsilon_h = h_x$.  }
\label{2d2_monotone_rates}
\end{table}

%%%

\begin{table}[htb] 
{\small
\begin{center}
\begin{tabular}{| c | c | c | c | c | c | c | c | c |}
		\hline
	& \multicolumn{4}{|c|}{Moment BC1} & \multicolumn{4}{|c|}{h*Moment BC1} \\ 
		\hline 
	 $h$ & $\ell^2$ Error & Order & $\ell^\infty$ Error & Order & $\ell^2$ Error & Order & $\ell^\infty$ Error & Order \\
		\hline
	1.57e-01 & 1.86e-02 &  & 3.87e-02 &  & 6.35e-03 &  & 2.02e-02 &  \\
		\hline
	7.44e-02 & 7.23e-03 & 1.26 & 2.22e-02 & 0.74 & 1.68e-03 & 1.78 & 8.61e-03 & 1.14 \\ 
		\hline
	3.63e-02 & 2.96e-03 & 1.24 & 1.27e-02 & 0.78 & 3.95e-04 & 2.01 & 2.81e-03 & 1.56 \\ 
		\hline
	2.40e-02 & 1.75e-03 & 1.27 & 8.75e-03 & 0.90 & 1.67e-04 & 2.08 & 1.36e-03 & 1.75 \\ 
		\hline
	1.79e-02 & 1.19e-03 & 1.31 & 6.57e-03 & 0.98 & 9.11e-05 & 2.08 & 7.99e-04 & 1.82 \\ 
		\hline
	1.19e-02 & 6.88e-04 & 1.34 & 4.29e-03 & 1.04 & 3.89e-05 & 2.08 & 3.71e-04 & 1.87 \\ 
		\hline
	8.89e-03 & 4.63e-04 & 1.36 & 3.11e-03 & 1.11 & 2.13e-05 & 2.07 & 2.13e-04 & 1.91 \\ 
		\hline
	5.92e-03 & 2.64e-04 & 1.38 & 1.95e-03 & 1.15 & 9.22e-06 & 2.06 & 9.67e-05 & 1.94 \\ 
		\hline
\end{tabular}
\end{center}
}
\caption{
Rates of convergence for Example 2 in two dimensions using the numerical moment and boundary condition \eqref{bc2a} 
with $\epsilon_h = h_x^2, \gamma_h = 4$ and $\epsilon_h = h^2, \gamma_h = 4 h_x$.  }
\label{2d2_momentBC1_rates}
\end{table}

%%%%%

\begin{table}[htb] 
{\small
\begin{center}
\begin{tabular}{| c | c | c | c | c | c | c | c | c |}
		\hline
	& \multicolumn{4}{|c|}{Moment BC2} & \multicolumn{4}{|c|}{h*Moment BC2} \\ 
		\hline 
	 $h$ & $\ell^2$ Error & Order & $\ell^\infty$ Error & Order & $\ell^2$ Error & Order & $\ell^\infty$ Error & Order \\
		\hline
	1.57e-01 & 4.48e-03 &  & 1.36e-02 &  & 4.26e-03 &  & 1.52e-02 &  \\
		\hline
	7.44e-02 & 2.29e-03 & 0.90 & 9.11e-03 & 0.54 & 1.16e-03 & 1.74 & 7.36e-03 & 0.97 \\ 
		\hline
	3.63e-02 & 8.55e-04 & 1.37 & 5.11e-03 & 0.80 & 2.72e-04 & 2.02 & 2.73e-03 & 1.38 \\ 
		\hline
	2.40e-02 & 4.58e-04 & 1.51 & 3.58e-03 & 0.86 & 1.13e-04 & 2.12 & 1.39e-03 & 1.62 \\
		\hline
	1.79e-02 & 2.90e-04 & 1.56 & 2.75e-03 & 0.91 & 6.06e-05 & 2.14 & 8.40e-04 & 1.73 \\
		\hline
	1.19e-02 & 1.50e-04 & 1.60 & 1.86e-03 & 0.95 & 2.51e-05 & 2.15 & 4.01e-04 & 1.80 \\ 
		\hline
	8.89e-03 & 9.37e-05 & 1.63 & 1.40e-03 & 0.99 & 1.35e-05 & 2.14 & 2.34e-04 & 1.86 \\
		\hline
	5.92e-03 & 4.77e-05 & 1.66 & 9.33e-04 & 0.99 & 5.70e-06 & 2.12 & 1.08e-04 & 1.90 \\ 
		\hline
\end{tabular}
\end{center}
}
\caption{
Rates of convergence for Example 2 in two dimensions using the numerical moment and boundary condition \eqref{bc2b} 
with $\epsilon_h = h_x^2, \gamma_h = 4$ and $\epsilon_h = h_x^2, \gamma_h = 4 h_x$.  }
\label{2d2_momentBC2_rates}
\end{table}

%%%%%%%%%%%%%%%%%%%%%%%
%%%%%%%%%%%%%%%%%%%%%%%

\subsubsection{Example 3:  Nonlinear Lipschitz operator with a lower-regularity solution} 
Consider the nonlinear problem \eqref{HJ} with 
\begin{align*}
	H[u] \equiv |u_x| + 2u_x & = f \qquad \text{in } \Omega , \\ 
	u & = g \qquad  \text{on } \partial \Omega , 
\end{align*}
where 
\[ f(x,y)= \begin{cases} 
      -1  & \text{if}\quad  x\leq 0.2 , \\
      3  &\text{if}\quad  x>0.2
   \end{cases}
\]
and $g$ is chosen such that $u(x,y) = | x - 0.2 |$. 
The results can be found in Tables~\ref{2d2_monotone_rates}, \ref{2d2_momentBC1_rates}, and \ref{2d2_momentBC2_rates}.  
As expected, the rates of convergence are decreased since the viscosity solution is not smooth.  
We can see similar accuracy results for the Lax-Friedrich's method and the proposed method with either choice for the 
auxiliary boundary condition when choosing $\gamma_h$ to be $\mathcal{O}(1)$ despite the lower rates of convergence 
for the numerical moment methods due to the smaller error on some of the coarse meshes.  
If we choose $\gamma_h = h_x$, the methods with a numerical moment stabilizer 
feature a more accurate coarse mesh approximation and the same linear 
rates of convergence as the Lax-Friedrich's method.   

\begin{table}[htb] 
{\small
\begin{center}
\begin{tabular}{| c | c | c | c | c |}
		\hline
	& \multicolumn{4}{|c|}{Lax-Friedrichs} \\ 
		\hline 
	 $h$ & $\ell^2$ Error & Order & $\ell^\infty$ Error & Order \\
		\hline
	1.57e-01 & 1.58e-01 &  & 2.65e-01 &  \\ 
		\hline
	7.44e-02 & 1.18e-01 & 0.39 & 1.72e-01 & 0.58 \\ 
		\hline
	3.63e-02 & 7.15e-02 & 0.70 & 9.44e-02 & 0.83 \\ 
		\hline
	2.40e-02 & 4.98e-02 & 0.87 & 6.33e-02 & 0.97 \\ 
		\hline
	1.79e-02 & 3.81e-02 & 0.92 & 4.74e-02 & 0.99 \\ 
		\hline
	1.19e-02 & 2.59e-02 & 0.94 & 3.15e-02 & 1.00 \\ 
		\hline
	8.89e-03 & 1.96e-02 & 0.96 & 2.35e-02 & 1.00 \\ 
		\hline
	5.92e-03 & 1.33e-02 & 0.96 & 1.57e-02 & 1.00 \\ 
		\hline
\end{tabular}
\end{center}
}
\caption{
Rates of convergence for Example 3 in two dimensions using the Lax-Friedrich's method 
with $\epsilon_h = 2 h_x$.  }
\label{2d3_monotone_rates}
\end{table}

%%%

\begin{table}[htb] 
{\small
\begin{center}
\begin{tabular}{| c | c | c | c | c | c | c | c | c |}
		\hline
	& \multicolumn{4}{|c|}{Moment BC1} & \multicolumn{4}{|c|}{h*Moment BC1} \\ 
		\hline 
	 $h$ & $\ell^2$ Error & Order & $\ell^\infty$ Error & Order & $\ell^2$ Error & Order & $\ell^\infty$ Error & Order \\
		\hline
	1.57e-01 & 1.99e-01 &  & 3.40e-01 &  & 1.50e-01 &  & 2.80e-01 &  \\ 
		\hline
	7.44e-02 & 1.17e-01 & 0.71 & 2.01e-01 & 0.70 & 4.47e-02 & 1.62 & 8.18e-02 & 1.65 \\ 
		\hline
	3.63e-02 & 3.60e-02 & 1.64 & 5.77e-02 & 1.74 & 2.19e-02 & 0.99 & 4.14e-02 & 0.95 \\ 
		\hline
	2.40e-02 & 1.90e-02 & 1.54 & 3.60e-02 & 1.14 & 1.41e-02 & 1.06 & 2.70e-02 & 1.03 \\ 
		\hline
	1.79e-02 & 1.59e-02 & 0.60 & 3.03e-02 & 0.59 & 1.04e-02 & 1.05 & 2.00e-02 & 1.03 \\ 
		\hline
	1.19e-02 & 1.42e-02 & 0.29 & 2.50e-02 & 0.47 & 6.81e-03 & 1.04 & 1.32e-02 & 1.01 \\
		\hline
	8.89e-03 & 1.14e-02 & 0.75 & 2.05e-02 & 0.69 & 5.06e-03 & 1.02 & 9.82e-03 & 1.02 \\ 
		\hline
	5.92e-03 & 7.95e-03 & 0.88 & 1.51e-02 & 0.74 & 3.34e-03 & 1.02 & 6.54e-03 & 1.00 \\ 
		\hline
\end{tabular}
\end{center}
}
\caption{
Rates of convergence for Example 3 in two dimensions using the numerical moment and boundary condition \eqref{bc2a} 
with $\epsilon_h = 2 h_x^2, \gamma_h = 1$ and $\epsilon_h = 2 h^2, \gamma_h = h_x$.  }
\label{2d3_momentBC1_rates}
\end{table}

%%%%%

\begin{table}[htb] 
{\small
\begin{center}
\begin{tabular}{| c | c | c | c | c | c | c | c | c |}
		\hline
	& \multicolumn{4}{|c|}{Moment BC2} & \multicolumn{4}{|c|}{h*Moment BC2} \\ 
		\hline 
	 $h$ & $\ell^2$ Error & Order & $\ell^\infty$ Error & Order & $\ell^2$ Error & Order & $\ell^\infty$ Error & Order \\
		\hline
	1.57e-01 & 1.89e-01 &  & 3.72e-01 &  & 1.25e-01 &  & 2.58e-01 &  \\
		\hline
	7.44e-02 & 7.57e-02 & 1.23 & 1.58e-01 & 1.15 & 5.02e-02 & 1.22 & 1.09e-01 & 1.15 \\ 
		\hline
	3.63e-02 & 3.55e-02 & 1.05 & 8.15e-02 & 0.92 & 2.26e-02 & 1.11 & 5.24e-02 & 1.02 \\ 
		\hline
	2.40e-02 & 2.76e-02 & 0.61 & 6.77e-02 & 0.45 & 1.44e-02 & 1.08 & 3.44e-02 & 1.01 \\ 
		\hline
	1.79e-02 & 2.25e-02 & 0.70 & 5.79e-02 & 0.53 & 1.06e-02 & 1.06 & 2.57e-02 & 1.01 \\
		\hline
	1.19e-02 & 1.55e-02 & 0.90 & 4.19e-02 & 0.79 & 6.90e-03 & 1.05 & 1.71e-02 & 0.99 \\ 
		\hline
	8.89e-03 & 1.20e-02 & 0.89 & 3.37e-02 & 0.75 & 5.11e-03 & 1.03 & 1.29e-02 & 0.99 \\ 
		\hline
	5.92e-03 & 8.69e-03 & 0.79 & 2.56e-02 & 0.68 & 3.37e-03 & 1.02 & 8.59e-03 & 0.99 \\ 
		\hline 
\end{tabular}
\end{center}
}
\caption{
Rates of convergence for Example 3 in two dimensions using the numerical moment and boundary condition \eqref{bc2b} 
with $\epsilon_h = 2 h_x^2, \gamma_h = 1$ and $\epsilon_h = 2 h_x^2, \gamma_h = h_x$.  }
\label{2d3_momentBC2_rates}
\end{table}

%%%%%%%%%%%%%%%%%%%%%%%%%%%%%%%%%%%%%%%%%%%%%%%%%%
%%%%%%%%%%%%%%%%%%%%%%%%%%%%%%%%%%%%%%%%%%%%%%%%%%
%%%%%%%%%%%%%%%%%%%%%%%%%%%%%%%%%%%%%%%%%%%%%%%%%%
%%%%%%%%%%%%%%%%%%%%%%%%%%%%%%%%%%%%%%%%%%%%%%%%%%

\clearpage 

\section{Conclusion} \label{conc_sec}

In this paper we examined the impact of adding a numerical moment stabilization term 
to a simple central difference approximation for 
linear reaction-convection-diffusion equations and stationary Hamilton-Jacobi equations.  
By introducing a numerical moment, we were able to form relatively accurate approximations on a coarse mesh 
even for the convection-dominated problem 
while also reliably approximating a discontinuity at the boundary  and solutions with corners.  
The schemes were able to break the first order accuracy barrier inherent to monotone schemes.  
Furthermore, the simple methods in this paper can also formally be extended to the discontinuous Galerkin setting 
to allow more flexible meshing and potentially increased accuracy 
as seen in \cite{DWDG}.  

The results in this paper have applications to 
approximating fully nonlinear second order elliptic problems 
and complements the FD methods proposed in \cite{FDhjb} and \cite{Kellie} that utilize numerical 
moments for approximating uniformly elliptic problems.  
By overcoming the need for monotonicity using the numerical moment, convergent narrow-stencil 
FD methods are available for a wide class of elliptic problems.  
This paper focused on consistency and $\ell^2$-stability analysis for linear problems as well as numerically testing 
the impact of the numerical moment as a complement to the paper \cite{HJfiltered} that focussed on 
$\ell^\infty$-stability analysis for nonlinear problems using a modified version of the scheme in this paper.  
A future direction is directly proving admissibility, stability, and convergence to the underlying viscosity solution 
for stationary HJ equations to rigorously extend the results in \cite{FDhjb} for fully nonlinear uniformly elliptic problems 
to fully nonlinear first order problems and degenerate elliptic problems.

%%%%%%%%%%%%%%%%%%%%%%%%%%%%%%%%%%%%%%%%%%%%%%%%%%
%%%%%%%%%%%%%%%%%%%%%%%%%%%%%%%%%%%%%%%%%%%%%%%%%%
%%%%%%%%%%%%%%%%%%%%%%%%%%%%%%%%%%%%%%%%%%%%%%%%%%
%%%%%%%%%%%%%%%%%%%%%%%%%%%%%%%%%%%%%%%%%%%%%%%%%%

\bibliographystyle{elsarticle-num}

\end{document}